\documentclass[12pt]{article}

\usepackage{amsmath,amssymb,amsthm}
\usepackage {graphicx}

\newtheorem{theorem}{Theorem}[section]
\newtheorem{lemma}[theorem]{Lemma}
\newtheorem{remark}[theorem]{Remark}

\newtheorem {proposition}[theorem]{Proposition}
\newtheorem {definition}[theorem]{Definition}
\newtheorem{problem}[theorem]{Problem}
\def \RR{\mathbb R}
\def \NN{\mathbb N}
\def\({\left(}  \def\){\right)}
\def\[{\left[}  \def\]{\right]}

\def \beq {\begin {equation}}
\def \eeq {\end{equation}}
\def \OL {\overline}

\def \W {\widetilde}

\begin {document}
\begin {center}
\Large{
Absolute Minima of Potentials\\ of a Certain Class of
Spherical Designs}

\bigskip

Sergiy Borodachov

\bigskip

\small {\it Department of Mathematics, Towson University, Towson, MD, 21252}
\end {center}
\large {

\begin {abstract}
We use linear programming techniques to find points of absolute minimum over the unit sphere $S^{d}$ in $\RR^{d+1}$ of the total potential of a point configuration $\omega_N\subset S^{d}$ which is a spherical $(2m-1)$-design contained in the union of some $m$ parallel hyperplanes. The interaction between points is described by the kernel $K({\bf x},{\bf y})=f(\left|{\bf x}-{\bf y}\right|^2)$, where $\left|\ \!\cdot\ \!\right|$ is the Euclidean norm in $\RR^{d+1}$. The potential function $f$ is assumed to have a convex derivative $f^{(2m-2)}$. Points of minimum do not depend on $f$ and are those and only those which form exactly $m$ distinct dot products with points of $\omega_N$. The proof of this theorem was presented at a workshop at ESI in January 2022. Using this result, we find sets of universal minima of certain six configurations on higher-dimensional spheres. 
\end {abstract}

{\it Keywords:} Spherical design, sharp configuration, extrema of a potential, universal extrema, Gegenbauer polynomials, absolutely (completely) monotone potential, Hermite interpolation.

{\it MSC 2020:} 33C45, 41A05, 52B11, 52C99, 33D90, 31B99.

\section {Introduction and setting of the problem}\label {intro}

Let $S^{d}:=\{(x_1,\ldots,x_{d+1})\in \RR^{d+1} : x_1^2+\ldots+x_{d+1}^2=1\}$ denote the unit sphere in $\RR^{d+1}$.
The problem of finding absolute minimum points for the total potential over $S^d$ of a point configuration (regular in some sense) was earlier studied in a number of works. In particular, papers \cite {Sto1975circle,Sto1975,NikRaf2011,NikRaf2013} focused on a regular $N$-gon on $S^1$ and a regular simplex, regular cross-polytope, and cube on $S^d$. Here we extend the results of these papers for the absolute minimum to a class of spherical designs which appears to contain quite a large number of other known configurations.
We use a general approach to characterizing the absolute minima of the potential by first finding the Hermite interpolating polynomial for the potential function and then applying properties of Gegenbauer polynomials. This is the classical Delsarte-Yudin (or linear programming) method, see the work by Cohn and Kumar \cite {CohKum2007}, book \cite {BorHarSafbook}, and references therein.

One important application of this problem is the maximal polarization problem on the sphere, which requires maximizing the absolute minimum over $S^{d}$ of the total potential of $N$ points lying on $S^{d}$. When the solution to the maximal polarization problem is expected to be a regular configuration, knowing the minimum value of its potential over the sphere would provide a sharp upper bound to prove. Reviews of known results on polarization can be found, for example, in book \cite [Chapter 14]{BorHarSafbook} with most recent work reviewed in \cite {BorMinMax}.

Let $f:[0,4]\to (-\infty,\infty]$ be a function finite and continuous on the interval $(0,4]$ with $\lim\limits_{t\to 0^+}f(t)=f(0)$. We call $f$ a {\it potential function}. Later, we will specify additional assumptions on its derivative. Points ${\bf x}$ and ${\bf y}$ on $S^{d}$ interact via the kernel of the form $K({\bf x},{\bf y})=f(\left|{\bf x}-{\bf y}\right|^2)$. Given a point configuration $\omega_N=\{{\bf x}_1,\ldots,{\bf x}_N\}\subset S^d$, we define {\it the $f$-potential of $\omega_N$} as
\begin {equation}\label {p_s}
p_f({\bf x},\omega_N):=\sum\limits_{i=1}^{N}f(\left|{\bf x}-{\bf x}_i\right|^2),\ \ \ {\bf x}\in S^d,
\end {equation}
and study the following problem.
\begin {problem}\label {P1}
{\rm Find the quantity
\begin {equation}\label {min}
P_f(\omega_N;S^d):=\min\limits_{{\bf x}\in S^d}p_f({\bf x},\omega_N),
\end {equation}
and determine all points ${\bf x}^\ast\in S^d$, for which the minimum is attained in \eqref {min}.
}
\end {problem}

Let $\omega_N^\ast$ be the configuration of the vertices of a regular $N$-gon inscribed in $S^1$. We will call dual for $\omega_N^\ast$ the regular $N$-gon formed by the midpoints of the arcs that join any two consecutive vertices of $\omega_N^\ast$. 
Let $\omega_{d+1}^\ast$ be the set of vertices of a regular simplex inscribed in $S^{d-1}$. The configuration $-\omega_{d+1}^\ast$ will be called its dual. 
Let $\omega_{2d}^\ast:=\{\pm {\bf e}_1,\ldots,\pm {\bf e}_d\}$, $d\geq 3$, where $\{{\bf e}_1,\ldots,{\bf e}_d\}$ is the standard basis in $\RR^d$, be the set of vertices of the regular cross-polytope inscribed in $S^{d-1}$, and $U_d$ be the set $\left\{\(\pm \frac {1}{\sqrt {d}},\ldots,\pm \frac {1}{\sqrt {d}}\)\right\}\subset \RR^d$ of $2^d$ vertices of the cube inscribed in $S^{d-1}$. The configurations $\omega_{2d}^\ast$ and $U_d$ will be called dual for each other. 
Recall that $\omega^\ast_N$ is a spherical $(N-1)$-design, $\omega_{d+1}^\ast$ is a $2$-design, while $\omega_{2d}^\ast$ and $U_d$ are both $3$-designs. 

{\it Early work.} Stolarsky \cite {Sto1975circle,Sto1975} found the absolute minima and maxima for the Riesz potential function $f(t)=-t^{-s/2}$ for a certain range of $s<0$ and the vertices of a regular polygon, regular simplex, regular cross-polytope, and a cube as well as absolute minimum for $f(t)=t^{-s/2}$, $s>0$, and the vertices of a regular cross-polytope and a cube. Later, Nikolov and Rafailov \cite {NikRaf2011,NikRaf2013} extended these results to all values of $s\neq 0$ for all four configurations $\omega^\ast_N$, $\omega^\ast_{d+1}$, $\omega^\ast_{2d}$, and $U_d$ and potential function
\begin {equation}\label {Rafailov}
f(t)=\begin {cases}\ \ (t+C)^{-s/2}, & s>0, \cr
-(t+C)^{-s/2}, & s<0,\cr
\end {cases} \ \ \ C\geq 0.
\end {equation}
One of the following three cases was proved to be always true:
the potential of the configuration achieves its absolute minimum over $S^d$ at the points of the dual configuration and its absolute maximum at the points of the configuration itself, or vice versa, or
the potential of the configuration is constant over $S^d$. Which case holds depends on whether the derivative of order one more than the spherical design strengh of the configuration is positive, negative, or zero on $(0,4]$.

The paper is structured as follows. Section \ref {PoSD} contains known facts and definitions used further in the paper while Section \ref {review} reviews known results on universal minima. In Section \ref {main}, we state the general results, Theorems \ref {2m-1} and \ref {2m-1w}, and in Section~\ref {immediate}, we give their immediate consequences. We prove Theorems \ref {2m-1} and \ref {2m-1w} in Section \ref {proofconst}. Section \ref {DDP} contains a known auxiliary statement on the distribution of dot products that a universal minimum point forms with points of a given stiff configuration. In Sections \ref {symm}, \ref {stiff}, and \ref {S6}, we find universal minima of the three pairs of mutually dual stiff configurations considered in this paper.

\section {Preliminaries}\label {PoSD}

In this section, we recall relevant definitions and known facts.
An infinitely differentiable function $f:(a,b)\to \RR$ is called {\it completely monotone} on $(a,b)$ if $(-1)^k f^{(k)}$ is non-negative on $(a,b)$ for all $k\geq 0$ and {\it strictly completely monotone} if $(-1)^k f^{(k)}$ is strictly positive on $(a,b)$ for all $k\geq 0$. Such are potential functions \eqref {Rafailov} for $s>0$ and $\W f(t)=e^{-at}$ for $a>0$. Function \eqref {Rafailov} for $C=0$ and $s>0$ defines the Riesz $s$-kernel $K_s({\bf x},{\bf y}):=\left|{\bf x}-{\bf y}\right|^{-s}$ and function $\W f$ the Gaussian kernel $K({\bf x},{\bf y})=e^{-a \left|{\bf x}-{\bf y}\right|^2}$. 
Furthermore, potential functions \eqref {Rafailov} for $-2<s<0$ and $f_{\log}(t)=\frac {1}{2}\ln \frac {1}{t}$ become strictly completely monotone on $(0,4)$ after adding an appropriate positive constant. Function \eqref {Rafailov} for $C=0$ and $s<0$ gives a rise to the Riesz kernel $K_s({\bf x},{\bf y}):=-\left|{\bf x}-{\bf y}\right|^{-s}$ and function $f_{\log}$ to the logarithmic kernel $K_{\log}({\bf x},{\bf y}):=\ln \frac {1}{\left|{\bf x}-{\bf y}\right|}$.

Completely monotone $f$ on $(0,4)$ (modulo an additive constant) have {\it convex} derivatives of even orders and {\it concave} derivatives of odd orders. Extended to $[0,4]$ by limits at the endpoints, they satisfy the assumptions of Theorems \ref {2m-1} and \ref {2m-1w}.

Denote 
$$
w_{d}(t):=\gamma_{d}(1-t^2)^{d/2-1}, 
$$
where the constant $\gamma_{d}$ is chosen so that $w_{d}$ is a probability density on $[-1,1]$. 
We recall that the {\it Gegenbauer (or ultraspherical) orthogonal polynomials} corresponding to the sphere $S^{d}$ in $\RR^{d+1}$ are terms of the sequence $\{P_n^{(d)}\}_{n=0}^{\infty}$ of univariate polynomials such that ${\rm deg}\ \! P_n^{(d)}=n$, $n\geq 0$, and 
$$
\int_{-1}^{1}P_i^{(d)}(t)P_j^{(d)}(t)w_{d}(t)\ \! dt =0,\ \ \ i\neq j,
$$
normalized so that $P_n^{(d)}(1)=1$, $n\geq 0$ (see \cite [Chapter 4]{Sze1975} or \cite [Chapter 5]{BorHarSafbook}). When $d=1$ or $d=2$, we obtain the classical Chebyshev polynomials of the first kind or the Legendre polynomials, respectively. When $d=3$, we have the Chebyshev polynomials of the second kind. It is a well-known fact that $P_n^{(d)}$ has exactly $n$ distinct zeros in $(-1,1)$ (see~\cite [Theorem 3.3.1]{Sze1975}), which we will denote by 
\begin {equation}\label {fund}
-1<\kappa_1^n<\ldots<\kappa_n^n<1.
\end {equation}
Since $P_n^{(d)}$ is even for even $n$ and odd for odd $n$, its set of roots is symmetric about $0$.

Let $\sigma_{d}$ denote the $d$-dimensional area measure on $S^d$ normalized to be a probability measure; that is, $\sigma_{d}\(\cdot\)=\frac {1}{\mathcal H_{d}(S^{d})}\mathcal H_d|_{S^{d}}\(\cdot\)$, where $\mathcal H_{d}$ is the $d$-dimensional Hausdorff measure in $\RR^{d+1}$.
Following the work by Delsarte, Goethals, and Seidel \cite {DelGoeSei1977}, we call a point configuration $\omega_N$ a {\it spherical $n$-design} if, for every polynomial $q$ on $\RR^{d+1}$ of degree at most~$n$,
\begin {equation}\label {design}
\frac {1}{N}\sum\limits_{i=1}^{N}q({\bf x}_i)=\int_{S^{d}}q({\bf x})\ \! d\sigma_{d}({\bf x}).
\end {equation}
The number $n$ is called the {\it strength} of the spherical design $\omega_N$.

Let $\mathbb P_n$ denote the space of all univariate polynomials of degree $\leq n$. An equivalent definition of a spherical $n$-design, $n\geq 1$, on $S^d$ is the fact \cite {DelGoeSei1977} that the equality  
\begin {equation}\label {constant}
\sum\limits_{i=1}^{N}p({\bf x}\cdot {\bf x}_i)=a_0(p)N,\ \ \ {\bf x}\in S^d,
\end {equation}
holds for every polynomial $p\in \mathbb P_n$. Here, $a_0(p)$ is the $0$-th Gegenbauer coefficient of $p$.


Let $m\in \NN$, ${\bf z}\in S^d$, and $-1\leq t_1<\ldots<t_m\leq 1$ be the distinct values of the dot product that ${\bf z}$ forms with points of the $n$-design $\omega_N$. Denote by $M_j$ the number of points in $\omega_N$ that form dot product $t_j$ with ${\bf z}$, $j=1,\ldots,m$.
Substituting ${\bf z}$ for ${\bf x}$ in \eqref {constant} we obtain a quadrature with nodes $t_j$ and weights $\beta_j:=M_j/N$ satisfying
\begin {equation}\label {cc}
a_0(p)=\int_{-1}^{1}p(t)w_d(t) \ \!dt=\sum\limits_{j=1}^{m}\beta_jp(t_j),\ \ \ p\in \mathbb P_n.
\end {equation}
Denote by $\{\varphi_1,\ldots,\varphi_m\}$ the fundamental set of polynomials for nodes $t_1,\ldots,t_m$; that is, $\varphi_k\in \mathbb P_{m-1}$, $\varphi_k(t_k)=1$, and $\varphi_k(t_i)=0$, $i\neq k$, $k=1,\ldots,m$. 

When $n\geq m-1$, the weights $\beta_j$ in \eqref {cc} are uniquely defined (see, e.g., \cite {IsaKel1965}): $\beta_j=a_0(\varphi_j)$, $j=1,\ldots,m$. To see this, simply let $p=\varphi_j$ in \eqref {cc}. Such a quadrature \eqref {cc} is called interpolatory. Thus, $M_j=a_0(\varphi_j)N$, $j=1,\ldots,m$; i.e., frequencies $M_j$ are determined by dot products $t_1,\ldots,t_m$ and the cardinality $N$ and do not depend on ${\bf z}$. Frequencies $M_j$ can also be found as solutions to a certain  Vandermonde linear system, see \cite [Theorem 2.1]{Boy1995}.

We next recall a well-known classification of quadratures with the highest algebraic degree of precision and dot product distributions for the corresponding classes of spherical designs (cf. \cite {DelGoeSei1977}). These quadratures are special cases of Gauss-Jacobi mechanical quadrature, see \cite[Theorem~3.4.1]{Sze1975}. The strict positivity of weights $\beta_j$ (called Christoffel numbers) in quadratures described below follows from \cite [Theorem 3.4.2]{Sze1975}.

{\it Case I}: $n=2m-1$. We will call $\omega_N$ $m$-stiff in this case, see Definition~\ref {D1}. Quadrature \eqref {cc} is now exact on $\mathbb P_{2m-1}$. 
The set $\{t_1,\ldots,t_m\}$ is the set of zeros of $P_m^{(d)}$. To see this, simply let $p(t)=\varphi_k(t)P_m^{(d)}(t)$ and use the equality in \eqref {cc}. 
Quadrature \eqref {cc} is called the Gauss-Gegenbauer quadrature.
\begin {remark}\label {m1}
{\rm 
From the above mentioned well-known facts, we immediately have the following. For every point ${\bf z}\in S^d$ that forms $m$ distinct dot products with points of a given $m$-stiff configuration, the values of dot products are zeros of the $m$-th Gegenbauer polynomial $P_m^{(d)}$ and their frequencies are independent of ${\bf z}$ and can be computed as $M_j=a_0(\varphi_j)N$.
}
\end {remark}

{\it Case II}: $n=2m-2$, $m\geq 2$, and ${\bf z}\in \omega_N$. Then quadrature \eqref {cc} is exact on $\mathbb P_{2m-2}$ with $t_m=1$ and $\beta_m=1/N$. The nodes $t_1,\ldots,t_{m-1}$ are zeros of the Jacobi orthogonal polynomial $J^{(\alpha,\beta)}_{m-1}$ of degree $m-1$, where $\alpha=d/2$ and $\beta=d/2-1$. To see this, let $p(t)=\varphi_k(t)J_{m-1}^{(\alpha,\beta)}(t)$, $k=1,\ldots,m-1$, in \eqref {cc} (which becomes $0$) and move the factor $1-t$ from $\varphi_k$ to the weight $w_d$ in the integral.

{\it Case III}: $n=2m-3$, $m\geq 2$, ${\bf z}\in \omega_N$, and $\omega_N$ is antipodal. Then quadrature \eqref {cc} is exact on $\mathbb P_{2m-3}$ with $t_1=-1$, $t_m=1$, and $\beta_1=\beta_m=1/N$. For $m\geq 3$, the nodes $t_2,\ldots,t_{m-1}$ are zeros of the Gegenbauer orthogonal polynomial $P^{(d+2)}_{m-2}$. To see this, let $p(t)=\varphi_k(t)P_{m-2}^{(d+2)}(t)$, $k=2,\ldots,m-1$, in \eqref {cc} and move factors $1-t$ and $1+t$ from $\varphi_k$ to the weight in the integral.

\begin {remark}\label {shar}
{\rm
If $\omega_N$ is a strongly $m$-sharp or $m$-sharp antipodal configuration on $S^d$ (see Definition \ref {sharp1}), then each point in $\omega_N$ forms the same set of $m$ dot products with other points of $\omega_N$ with the same frequencies (which can be computed as described above). This set of dot products is the set of zeros of the Jacobi polynomial $J_m^{(\alpha,\beta)}$, where $\alpha=d/2$ and $\beta=d/2-1$, when $\omega_N$ is strongly $m$-sharp and is the set of zeros of the Gegenbauer polynomial $P_{m-1}^{(d+2)}$ together with $-1$ when $\omega_N$ is $m$-sharp antipodal.
}
\end {remark}

\section {Review of known results on universal extrema}\label {review}

A number of recent works were devoted to finding extrema of general potentials of regular configurations.
Earlier results for power-law potentials were reviewed in Section \ref {intro}.
Among kernels of the form $K({\bf x},{\bf y})=f(\left|{\bf x}-{\bf y}\right|^2)$, of greatest interest here are the ones with a completely monotone potential function $f$ on $(0,4)$ or with $f$ that becomes completely monotone after adding a constant.

\begin {definition}
{\rm
Given a configuration $\omega_N\subset S^d$, a point ${\bf x}^\ast\in S^d$ is called a  universal minimum point (of the potential) of $\omega_N$ if ${\bf x}^\ast$ attains the absolute minimum in \eqref {min} for every potential function $f$ completely monotone on $(0,4]$.
}
\end {definition}
For universal minimality results obtained here, we will require only that the derivative of $f$ of a certain even order be convex or the derivative of a certain odd order be concave.
Hardin, Kendall, and Saff \cite {HarKenSaf2013} proved that absolute minima of the potential of $\omega_N^\ast$ on $S^1$ 
with respect to a decreasing and convex function of the geodesic distance are attained at points of the dual regular $N$-gon. This partially solved the problem about universal minima of $\omega^\ast_N$. The set of universal minima of regular simplex $\omega_{d+1}^\ast$ was shown in \cite {Borsimplex} to be $-\omega_{d+1}^\ast$. 

During talk \cite {Bor2022talk} given in January 2022 at ESI (available in the official YouTube account of ESI) the author of this paper proved the following theorem: given an $m$-stiff configuration $\omega_N$ (see Definition \ref {D1}), its universal minima are those and only those points on $S^d$ that form $m$ distinct dot products with points of $\omega_N$. We state this result in Theorem \ref {2m-1}\footnote {Theorem \ref {2m-1} was earlier announced in talk \cite {BorAIMtalk} in May 2021.}. In view of Remark \ref {m1} from the preliminaries, one can write
\begin {equation}\label {e1}
\min\limits_{{\bf x}\in S^d}p_f({\bf x},\omega_N)=\sum\limits_{i=1}^{m}M_if(2-2t_i),
\end {equation}
where $t_1,\ldots,t_m$ are zeros of $P_m^{(d)}$, $\{\varphi_1,\ldots,\varphi_m\}$ is the correponding fundamental set of polynomials, and  $M_i=a_0(\varphi_i)N$, $i=\OL{1,m}$.

During talk \cite {Bor2022talk}, we mentioned regular $2m$-gon on $S^1$, cube and cross-polytope on $S^{d-1}$, and the $24$-cell on $S^3$ as examples of stiff configurations.
Since the duality of the cube and cross-polytope can be easily seen as well as the duality of certain two isometric copies of the $24$-cell, this theorem immediately implies the following. The set of universal minima of a regular $2m$-gon on $S^1$ is the set of midpoints of the arcs joining its neighboring points. The set of universal minima of a regular cross-polytope is the dual cube and the one for a cube is the dual cross-polytope. 
The set of universal minima of the $24$-cell is the dual $24$-cell. We state these corollaries in Section \ref {immediate}.

The above mentioned theorem on stiff configurations is a straightforward consequence of Lemma 3.5 in \cite {BorMinMax} (posted to ArXiv in March 2022). We state it here as Lemma~\ref {gen}. Lemma 3.5 in \cite{BorMinMax} has three other similar immediate consequences. All four follow if one just uses the well-known error formula for Hermite interpolation (see, e.g., \cite [Section 2.2]{IsaKel1965} or \cite[Theorem~4.5]{Epp2013}) or its ``slight" modifications to include interpolation at the endpoint(s) (see, e.g., \cite [Theorem 3.5.1] {Dav1975}). One of those consequences\footnote {Mentioned in abstract \cite [p. 78]{BorBatumi} accepted on July 1, 2022.} is that, for any strongly sharp configuration $\omega_N$ (see Definition~\ref {sharp1}), the set of its universal minima is the set of antipods of $\omega_N$. We state it in Theorem~\ref {2m-1w}. Therefore, the set of universal minima of a regular $(2m+1)$-gon is the set of midpoints of the arcs joining neighboring vertices and the sets of universal minima of the Schl\"affi configuration on $S^5$ and of the McLaughlin configuration on $S^{21}$ are at the antipods of points of these configurations (universal minima of regular simplex were mentioned earlier).

The remaining two consequences of \cite [Lemma 3.5]{BorMinMax} are that universal maxima of any strongly sharp or sharp antipodal configuration (see Definition~\ref {sharp1}) are (only) at the points of the configuration itself, see \cite {BorMinMax}\footnote{The proof of this was also given earlier in talk \cite {Bor2022talk}.}. In particular,
universal maxima of a regular $N$-gon, regular simplex, regular cross-polytope, regular icosahedron, and $E_8$ and Leech lattices are points of the configuration itself.

Dot products $t_1,\ldots,t_m$ between distinct points in a strongly $m$-sharp or $m$-sharp antipodal configuration are well-known (see Remark \ref {shar}). The corresponding frequencies $M_i$ (independent of the point in $\omega_N$) can be computed as explained in Section \ref {PoSD}. So one can write
\begin {equation}\label {e2}
\min\limits_{{\bf x}\in S^d}p_f({\bf x},\omega_N)=f(4)+\sum\limits_{i=1}^{m}M_if(2+2t_i),
\end {equation}
for the minimum value for a strongly $m$-sharp configuration and
\begin {equation}\label {e3}
\max\limits_{{\bf x}\in S^d}p_f({\bf x},\omega_N)=f(0)+\sum\limits_{i=1}^{m}M_if(2-2t_i),
\end {equation}
for the maximum value for a strongly $m$-sharp or an $m$-sharp antipodal one, where in the case of a sharp antipodal configuration, we have $t_1=-1$ and $M_1=1$.

The work described above stimulated further research. The universal maxima of the $600$-cell on $S^3$ (which is not sharp) were shown by Boyvalenkov, Dragnev, Hardin, Saff, and Stoyanova~\cite {BoyDraHarSafSto600cell} to be at points of the $600$-cell itself (on ArXiv since July 15, 2022). 
A number of new stiff configurations (which are also sharp) were found by Boyvalenkov, Dragnev, Hardin, Saff, and Stoyanova \cite {BoyDraHarSafStosharpantipodal} (on ArXiv since October 2022). The result of \cite {Bor2022talk} implies that their universal minima are at points of the corresponding dual configurations (see Definition~\ref {D1}). Paper \cite {BoyDraHarSafStosharpantipodal} further elaborates on each of these dual configurations (as well as on antipods of strongly sharp configurations on $S^{21}$ mentioned above).


For stiff or strongly sharp configurations considered in \cite {BoyDraHarSafStosharpantipodal}, values of expressions in \eqref {e1} or \eqref {e2} are computed in \cite [Table 2] {BoyDraHarSafStosharpantipodal}. 
Results of \cite [Theorem 1.4]{BoyDraHarSafStosharpantipodal} (see also \cite {BoyDraHarSafSto600cell}) provide lower and upper bounds for the potential of any configuration $\omega_N$ on $S^d$ as soon as $\omega_N$ is at least a $1$-design. When $\omega_N$ is stiff, strongly sharp, or sharp antipodal, these bounds become sharp (one essentially gets equations \eqref {e1}--\eqref {e3} mentioned above\footnote {In \cite[Theorem 1.4] {BoyDraHarSafStosharpantipodal} item (i) case $\epsilon=1$, dot products $\alpha_j$ are zeros of $(1+t)^\epsilon P_k^{(0,\epsilon)}(-t)$.}). Proofs of Lemma~3.5 from \cite {BorMinMax} and of Theorem~1.4 from \cite {BoyDraHarSafStosharpantipodal} both use equivalent definition \eqref {constant} of a spherical design.

Some important regular configurations are not stiff or strongly sharp. The author \cite {Borsymmetric} (on ArXiv since October 9, 2022) found and characterized all universal minima for the icosahedron (minima are only at vertices of the dual dodecahedron), dodecahedron (dual icosahedron), and for the $E_8$ lattice on $S^7$ ($2_{41}$ polytope). This was done by proving a general theorem for $(2m-3)$-designs with vanishing $(2m-1)$-th and $2m$-th moments (the ``skip one add two" case). It was explicitly stated in talk \cite {BorBatumi}\footnote{The talk was given on August 29, 2022.} and its proof was earlier illustrated using the example of icosahedron in talk~\cite {Bor2022talk}\footnote{Dodecahedron also fits the argument for icosahedron from talk \cite {Bor2022talk}.}.

One universal minimum and the universal minimum value of the potential of the Leech lattice on $S^{23}$ was found by Boyvalenkov, Dragnev, Hardin, Saff, and Stoyanova \cite {BoyDraHarSafStosharpantipodal} (on ArXiv since October 31, 2022) who also proved a general statement for the ``skip one add two" case. Their theorem resorts to the Gaussian-type quadrature for the corresponding polynomial space. It provides a different (to a certain extent) proof for universal minima of icosahedron, dodecahedron, and (without characterization of minima) $E_8$ lattice.

Critical points of the total potential of finite configurations of charges were also analysed (see \cite {Bil2015,GioKhi2020} and references therein). This work is related to the known Maxwell's conjecture.

We remark that the problem about universal minima of the $600$-cell remains open as well as the problem about universal maxima of stiff configurations that are not sharp antipodal (e.g., cube).

In addition to Theorems \ref {2m-1} and \ref {2m-1w}, we find in this paper all universal minima of six stiff configurations, for which universal minima have not been characterized before, see Sections \ref {symm}, \ref {stiff}, and \ref {S6}. 

Table \ref {tableA} shows the current state of the art in finding universal extrema for the  classes of configurations of interest and for the $600$-cell. Table \ref {table1} summarizes results of our work in \cite {Borsymmetric,BorMinMax,Bor2022talk} and in this paper for particular configurations. Last three lines of Table \ref {table2} contain the six stiff configurations, for which we find here universal minima (the first three lines follow immediately from Theorem~\ref {2m-1}).

\begin {table}
\centering
\begin{tabular}{|c|c|c|}
\hline
 {\bf Class of configurations} & {\bf Universal maxima}  &  {\bf Universal minima } \\
\hline
  {Strongly sharp}  & Solved \cite {BorMinMax} & Theorem \ref {2m-1w}, announced \cite {BorBatumi}, see also \cite {BoyDraHarSafStosharpantipodal}  \\
\hline
 {Sharp antipodal}, {non-stiff} & Solved \cite {BorMinMax} & Solved \cite {Borsymmetric,BoyDraHarSafStosharpantipodal} for a subclass \\
\hline
 {Sharp antipodal and stiff} & Solved \cite {BorMinMax} & Solved \cite{Bor2022talk}, see Theorem \ref {2m-1} \\
\hline

Stiff, not sharp antipodal & Open & Solved \cite {Bor2022talk}, see Theorem \ref {2m-1}\\

\hline
 $600$-cell & Solved \cite {BoyDraHarSafSto600cell} & Open \\
\hline

\end {tabular}
\caption {Current progress on the solution of the problem about universal extrema of general classes of spherical configurations of interest and of the $600$-cell.}\label {tableA}
\end {table}

When proving our main results, we use the Delsarte method (the Delsarte-Yudin bound, linear programming method, or polynomial method), see the work by Delsarte, Goethals, and Seidel \cite {DelGoeSei1977}. This approach has been extensively used, in particular, to find best-packing, kissing number, and minimal energy, as well as to estimate the number of points in a spherical design (see, \cite {Lev1979,Lev1992,Lev1998,BoyDanKaz2001,Mus2003,CohKum2007,BorHarSafbook} and references therein).

\section {General results}\label {main}

We first state the solution to Problem \ref {P1} for a class of spherical designs, which we call here stiff configurations.
Throughout the rest of the paper, ${\bf x}_1,\ldots,{\bf x}_N$ will denote the points in a configuration $\omega_N$ on $S^{d}$ whose total potential we are minimizing.
For a given $m\in \NN$ and a given configuration $\omega_N\subset S^{d}$,
denote by $\mathcal D_m(\omega_N)$ the set of all points ${\bf z}\in S^{d}$ for which the set of dot products 
$$
D({\bf z},\omega_N):=\{{\bf z}\cdot {\bf x}_i : i=1,\ldots,N\}
$$ 
has at most $m$ distinct elements. The set $D({\bf z},\omega_N)$ is allowed to contain $1$, which occurs when ${\bf z}\in \omega_N$. If $\mathcal D_m(\omega_N)\neq \emptyset$, then $\omega_N$ is contained in a union of at most $m$ parallel hyperplanes.

\begin {definition}\label {D1}
{\rm
We call a point configuration $\omega_N\subset S^{d}$ {\it $m$-stiff}, $d,m\geq 1$, if $\omega_N$ is a spherical $(2m-1)$-design and the set $\mathcal D_m(\omega_N)$ is non-empty. The set $\mathcal D_m(\omega_N)$ of a given $m$-stiff configuration $\omega_N$ is called {\it the dual configuration} for~$\omega_N$.
}
\end {definition}
\begin {remark}\label {m}
{\rm
If $\omega_N$ is an $m$-stiff configuration on $S^d$, $d,m\geq 1$, then every point in $\mathcal D_m(\omega_N)$ forms exactly $m$ distinct dot products with points of $\omega_N$. If it did not, then $\omega_N$ would be contained in $k<m$ parallel hyperplanes ${\bf x}\cdot {\bf v}=\alpha_i$, $i=1,\ldots,k$, and the polynomial $u({\bf x}):=({\bf x}\cdot {\bf v}-\alpha_1)^2\cdots ({\bf x}\cdot {\bf v}-\alpha_k)^2$ of degree $2k<2m-1$ would vanish on $\omega_N$ while having a positive integral over $S^d$.
}
\end {remark}

In view of this remark, the dual configuration for an $m$-stiff $\omega_N$ can be equivalently defined as the set of all points on $S^d$ that form exactly $m$ distinct dot products with points of $\omega_N$.

There are many examples of stiff configurations other than the ones mentioned when reviewing known results in Section \ref {intro}. First, there are two interesting examples on $S^{d-1}$ for $d\geq 3$ odd (both are $2$-stiff), two more on $S^5$, and two more on $S^6$ (all are $3$-stiff). We describe these configurations in Sections~\ref {symm}, \ref  {stiff}, and \ref {S6}. Furthermore, paper \cite {BoyDraHarSafStosharpantipodal} analyses one example on $S^4$ and a number of examples on higher-dimensional spheres.

We have the following result.
It was prsented in talk \cite {Bor2022talk}\footnote {Given in January 2022.} with the proof. 
\begin {theorem}\label {2m-1}
Let $m\geq 1$, $d\geq 1$, $f:[0,4]\to (-\infty,\infty]$ be a function continuous on $(0,4]$ with $f(0)=\lim\limits_{t\to 0^+}f(t)$, differentiable on $(0,4)$ with a convex derivative $f^{(2m-2)}$ on $(0,4)$. If $\omega_N=\{{\bf x}_1,\ldots,{\bf x}_N\}$ is an $m$-stiff configuration on the sphere $S^{d}$, then the potential
$$
p_f({\bf x},\omega_N)=\sum\limits_{i=1}^{N}f(\left|{\bf x}-{\bf x}_i\right|^2),\ \ \ {\bf x}\in S^{d}, 
$$
attains its absolute minimum over $S^{d}$ at every point of the set $\mathcal D_m(\omega_N)$. 

If, in addition, $f^{(2m-2)}$ is strictly convex on $(0,4)$, then $\mathcal D_m(\omega_N)$ contains all points of absolute minimum of the potential $p_f(\cdot,\omega_N)$ on $S^{d}$.
\end {theorem}

Immediate consequences of Theorem \ref {2m-1} are stated in Section \ref {immediate}.
No point from an $m$-stiff configuration $\omega_N$ (or $-\omega_N$) can be 
in $\mathcal D_m(\omega_N)$, see Remark \ref {m1} or Lemma \ref {w2m-1}. Therefore, a point from $\omega_N$ in Theorem \ref {2m-1} cannot be a minimizer of the potential of $\omega_N$ when $f^{(2m-2)}$ is strictly convex. 

The ``slightly" less general case of Theorem \ref {2m-1} when $f$ has a positive derivative $f^{(2m)}$ on $(0,4]$ is an important special case of our earlier result \cite[Lemma~3.5]{BorMinMax} restated here as Lemma \ref {gen}. To see this, recall the following well-known fact from a course in numerical analysis \cite [Section 2.2]{IsaKel1965} or \cite [Theorem 4.5]{Epp2013}. Let $h$ denote the Hermite interpolating polynomial for $g(t)=f(2-2t)$ at the set of nodes $-1< t_1<\ldots<t_m< 1$, which is chosen to be the set of dot products that a point from $\mathcal D_m(\omega_N)$ forms with points from $\omega_N$. Then ${\rm deg}\ \! h\leq 2m-1$ and
\begin {equation}\label {H}
g(t)-h(t)=\frac{g^{(2m)}(c)}{(2m)!}(t-t_1)^2\cdots(t-t_m)^2\geq 0,\ \ \ t\in [-1,1),
\end {equation}
where $c=c(t)\in (-1,1)$. If $f(1)<\infty$, then \eqref {H} also holds for $t=1$.
The ``full" case of Theorem \ref {2m-1} follows on applying Lemma~\ref {Interp} instead of \eqref {H}, see details in Section \ref {proofconst}.

We consider next another class of regular configurations on the sphere, for which universal minima can be found. If the set $\mathcal D_{m}(\omega_N)$ contains a point from $-\omega_N$, then it is enough to require that $\omega_N$ be a $(2m-2)$-design. These assumptions are satisfied by any sharp configuration, which is an even order design. For the definition of a sharp configuration, see~\cite {CohKum2007}.
\begin {definition}\label {sharp1}
{\rm
We call a point configuration $\omega_N\subset S^{d}$ {\it $m$-sharp}, $m\geq 1$, if $\omega_N$ is a spherical $(2m-1)$-design and there are $m$ distinct values of the dot product that can occur between any two distinct vectors from $\omega_N$. We call an $m$-sharp configuration {\it strongly $m$-sharp} if it is a spherical $2m$-design.
}
\end {definition}

Examples of strongly sharp configurations are vertices of a regular $(2m+1)$-gon inscribed in $S^1$, $m\in \NN$, the set of vertices of a regular $d$-simplex inscribed in $S^{d-1}$ ($m=1$), the Schl\"affi configuration of $N=27$ points on $S^5$ ($m=2$), and the McLaughlin configuration of $N=275$ points on $S^{21}$ ($m=2$). There are no strongly $m$-sharp configurations on $S^d$ for $d\geq 2$ and $m\geq 3$. This follows from the results of \cite {DelGoeSei1977} and \cite {BanDam1979}.

The classes of stiff and sharp configurations overlap. For example, sets of vertices of a regular cross-polytope on $S^{d-1}$ and of a regular $2m$-gon inscribed in $S^1$ are both stiff and sharp as well as the kissing configuration $\OL\omega_{56}$ on $S^6$ and some examples mentioned in \cite [Table 2]{BoyDraHarSafStosharpantipodal}. However, there are stiff configurations, which are not sharp (e.g., the cube in $\RR^d$ for $d\geq 3$ or the $24$-cell on $S^3$). There are also sharp configurations, which are not stiff, e.g., the regular $(2m+1)$-gon inscribed in $S^1$, the regular icosahedron on~$S^2$, and the minimal vectors of $E_8$ lattice on $S^7$. In addition, classes of strongly sharp and stiff configurations are disjoint (see Proposition~\ref {sst}).
%
%
%

The following result\footnote {Mentioned in the abstract of~\cite [p. 78]{BorBatumi}. } follows from \cite [Lemma 3.5]{BorMinMax} in the same way as Theorem \ref {2m-1}. 

\begin {theorem}\label {2m-1w}
Let $m\geq 1$, $d\geq 1$, $f:[0,4]\to (-\infty,\infty]$ be a function continuous on $(0,4]$ with $f(0)=\lim\limits_{t\to 0^+}f(t)$, differentiable on $(0,4)$ with a concave derivative $f^{(2m-1)}$ on $(0,4)$. If $\omega_N=\{{\bf x}_1,\ldots,{\bf x}_N\}\subset S^{d}$ is a strongly $m$-sharp configuration, then the potential
$$
p_f({\bf x},\omega_N)=\sum\limits_{i=1}^{N}f(\left|{\bf x}-{\bf x}_i\right|^2),\ \ \ {\bf x}\in S^{d}, 
$$
attains its absolute minimum over $S^{d}$ at every point of  $-\omega_N$.

If the concavity of $f^{(2m-1)}$ on $(0,4)$ is strict, then $-\omega_N$ contains all points of absolute minimum of the potential $p_f(\cdot,\omega_N)$ over $S^{d}$.
\end {theorem}

%
%
%
%
%
%

\begin {table}
\centering
\small{\begin{tabular}{ |c|c|c|c|c|c|c| } 
 \hline
Name & $N$ & Sphere & Design & Sharpness & Stiffness & We  \\ 
& & & strength & & & find \\
\hline
Regular $(2m+1)$-gon & $2m+1$ & $S^1$ & $2m$ & strongly $m$-sharp & $-$ & both\\
Regular $2m$-gon  & $2m$ & $S^1$ & $2m-1$ & $m$-sharp & $m$-stiff & both \\
Regular $d$-simplex & $d+1$ & $S^{d-1}$ & $2$ & strongly $1$-sharp & $-$ & both \\ 
Regular cross-polytope ($d\geq 3$)& $2d$ & $S^{d-1}$ & $3$ & $2$-sharp & $2$-stiff & both\\
Cube ($d\geq 3$) & $2^d$ & $S^{d-1}$ & $3$ & $-$ & $2$-stiff & min\\
Symmetrized $(2\ell+1)$-simplex & $4\ell+4$ & $S^{2\ell}$ & $3$ & $-$ & $2$-stiff & min\\
Dual symm.$(2\ell\!+\!1)$-simpl., $\ell\!\geq\! 2$ & $\binom{2\ell+2}{\ell+1}$ & $S^{2\ell}$ & $3$ & $-$ & $2$-stiff & min \\
Icosahedron & $12$ & $S^2$ & $5$ & $3$-sharp & $-$ & both\\ 
Dodecahedron & $20$ & $S^2$ & $5$ & $-$ & $-$ & min \\
$24$-cell & $24$ & $S^3$ & $5$ & $-$ & $3$-stiff & min \\
Schl\" affi & $27$ & $S^5$ & $4$  & strongly $2$-sharp & $-$ & both\\
Symmetrized Schl\"affi & $54$ & $S^5$ & $5$ & $-$ & $3$-stiff & min\\
Roots of $E_6$ lattice & $72$ & $S^5$ & $5$ & $-$ & $3$-stiff & min \\
Kissing & $56$ & $S^6$ & $5$ & $3$-sharp & $3$-stiff & both\\
$E_7$ lattice & $126$ & $S^6$ & $5$ & $-$ & $3$-stiff & min\\
$E_8$ lattice & $240$ & $S^7$ & $7$ & $4$-sharp & $-$ & both \\
McLaughlin & $275$ & $S^{21}$ & $4$ & strongly $2$-sharp & $-$ & both\\
Leech lattice & $196560$ & $S^{23}$ & $11$ & $6$-sharp & $-$ & max \\
\hline

\end{tabular}
}\caption{Summary of our results in \cite {Borsymmetric,BorMinMax,Bor2022talk} and in this paper giving all universal extrema of certain configurations (see the last column). Entry with a dash indicates that the property does not hold. We remark that all universal maxima of the $600$-cell were found in \cite {BoyDraHarSafSto600cell} and one universal minimum was found in \cite{BoyDraHarSafStosharpantipodal} for the kissing $56$-point configuration on $S^6$ and for Leech lattice on $S^{23}$. Paper \cite {BoyDraHarSafStosharpantipodal} studies universal minima for a number of stiff configurations (which are also sharp) not listed in this table. 
}\label {table1}. 
\end {table}

\begin {table}
\centering
\begin{tabular}{|c|c|c|}
\hline
{\bf Sphere} & \multicolumn {2}{c|}{\bf Pairs of mutually dual stiff configurations} \\
\hline
$S^1$ & regular $2m$-gon & regular $2m$-gon rotated by $\frac{\pi}{2m}$ rad\\
\hline
$S^{d-1}$, $d\geq 3$ & regular cross-polytope $\omega^\ast_{2d}$ & cube $U_d$ \\
\hline
$S^3$ & $24$-cell $\OL\omega_{24}$ & $24$-cell $\W\omega_{24}$\\
\hline
$S^{d-1}$, $d\geq 3$ odd & symmetrized regular $d$-simplex $\OL\omega_{2d+2}$ & $\OL\omega_{\kappa(d)}$ defined by \eqref {SymSimp}\\
\hline
$S^5$ & symmetrized Schl\"affi configuration $\OL\omega_{54}$ & root system of $E_6$-lattice $\OL\omega_{72}$ \\
\hline
$S^6$ & kissing configuration $\OL\omega_{56}$ & minimal vectors of $E_7$-lattice $\OL\omega_{126}$ \\
\hline
\end {tabular}
\caption {Pairs of mutually dual stiff configurations. Universal minima for the first three lines are immediate consequences of Theorem \ref {2m-1} proved in talk \cite {Bor2022talk}. Universal minima of configurations in the last three lines are obtained in Sections \ref {symm}, \ref {stiff}, and \ref {S6}.}\label {table2}
\end {table}

Theorem \ref {2m-1w} when $g(t)=f(2-2t)$ has a positive derivative $g^{(2m+1)}$ on $[-1,1)$ can be derived from \cite[Lemma 3.5]{BorMinMax} by recalling the following  fact from a numerical analysis course \cite [Theorem 3.5.1]{Dav1975}. If $h$ is the polynomial of degree at most $2m$ interpolating $g$ to the second order at the nodes $-1< t_1<\ldots<t_m<1$ and to the first order at $t_0=-1$, then
\begin {equation}\label {H1}
g(t)-h(t)=\frac{g^{(2m+1)}(c)}{(2m+1)!}(t-t_1)^2\cdots(t-t_m)^2(t+1)\geq 0,\ \ \ t\in [-1,1),
\end {equation}
where $c=c(t)\in (-1,1)$. Numbers $t_1,\ldots,t_m$ are chosen to be opposite to dot products between distinct points in $\omega_N$. If $f(1)<\infty$, then \eqref {H1} also holds for $t=1$.
The ``full" case of Theorem \ref {2m-1w} follows on applying Lemma~\ref {Interp} instead of \eqref {H1}, see details in Section \ref {proofconst}. 


Note that in the proofs of Theorems \ref {2m-1} and \ref {2m-1w}, we do not need to show the positive definiteness of the interpolating polynomial when utilizing the linear programming approach of \cite {DelGoeSei1977,CohKum2007}.

Other main results of this paper are Theorems \ref {symmetrizedsimplex}, \ref {72_}, and \ref {56_}.

\section {Immediate consequences of Theorems \ref {2m-1} and \ref {2m-1w}}\label {immediate}

This section discusses universal minima given by Theorem \ref {2m-1} for regular $2m$-gon, regular cross-polytope, cube, and $24$-cell (finding the dual for each configuration is elementary) and universal minima given by Theorem~\ref {2m-1w} for a regular $(2m+1)$-gon, regular simplex, and Schl\"affi and McLaughlin configurations.
For power potentials \eqref {Rafailov}, the first five lines of Table \ref {table4} were proved in \cite {Sto1975circle,Sto1975,NikRaf2011,NikRaf2013}. For regular simplex, universal minima were found in \cite {Borsimplex} (here we provide an alternative proof).

The following auxiliary statements give dual configurations.
Let $\W\omega_N$ denote the set of midpoints of arcs on $S^1$ that join the neighboring points in the regular $N$-gon $\omega^\ast_N$ inscribed in $S^1$.
\begin {lemma}\label {Wo}
For $N=2m$, $m\geq 1$, we have $\mathcal D_m(\omega_N^\ast)=\W\omega_N$.
\end {lemma}
\begin {proof}
Every point from $\W\omega_{N}$, $N=2m$, forms $m$ distinct dot products with points of $\omega_N^\ast$, while every point from $\omega_N^\ast$ itself forms $m+1$ dot products and every other point on $S^1$ forms $2m$ dot products with points of $\omega_N^\ast$.
\end {proof}

\begin {lemma}\label {omegacube}
For every $d\geq 2$, we have $\mathcal D_2(\omega_{2d}^\ast)=U_d$.
\end {lemma}
\begin {proof}
The inclusion $U_d\subset \mathcal D_2(\omega_{2d}^\ast)$ is immediate. Let $x_i$ be a non-zero coordinate of ${\bf x}:=(x_1,\ldots,x_d)\in \mathcal D_2(\omega_{2d}^\ast)$. Then ${\bf x}\cdot (\pm {\bf e}_i)=\pm x_i$. Every other coordinate $x_j$ must have the same absolute value as $x_i$ (otherwise $\pm x_j$ would give extra dot product(s)). Then $\left|x_i\right|=\frac {1}{\sqrt{d}}$, $i=\OL{1,d}$; that is, ${\bf x}\in U_d$. 
\end {proof}

\begin {lemma}\label {cubeomega}
For every $d\geq 2$, we have $\mathcal D_2(U_d)=\omega_{2d}^\ast$.
\end {lemma}
\begin {proof}
The inclusion $\omega_{2d}^\ast\subset \mathcal D_2(U_d)$ is immediate. If ${\bf x}=(x_1,\ldots,x_d)\in \mathcal D_2(U_d)$ had at least two non-zero coordinates, say $x_1$ and $x_2$, then dot product ${\bf x}\cdot \(\pm \frac {1}{\sqrt{d}},\pm \frac {1}{\sqrt{d}},\frac {1}{\sqrt{d}},\ldots,\frac {1}{\sqrt{d}}\)=\frac {\pm x_1\pm x_2}{\sqrt{d}}+\frac {x_3}{\sqrt{d}}+\ldots+\frac {x_d}{\sqrt{d}}$ would have at least three distinct values. Therefore, only one coordinate of ${\bf x}$ is non-zero; i.e., ${\bf x}\in \omega_{2d}^\ast$.
\end {proof}

The $24$-cell $\OL\omega_{24}\subset S^3$ is the set of $N=24$ points with two coordinates equal $0$ and each of the remaining two coordinates being $1/\sqrt {2}$ or $-1/\sqrt {2}$. The $24$-cell is known to be a $5$-design. At the same time, every vector $\pm{\bf e}_i$, $i=\OL {1,4}$, where $\{{\bf e}_1,{\bf e}_2,{\bf e}_3,{\bf e}_4\}$ is the standard basis in $\RR^4$, forms only $3$ distinct dot products with vectors of the $24$-cell, which are $\pm 1/\sqrt {2}$ and $0$. Furthermore, each of the $16$ vectors of the form $\(\pm \frac {1}{2},\pm \frac {1}{2},\pm \frac {1}{2},\pm \frac {1}{2}\)$ forms these three values of the dot product with vectors from $\OL \omega_{24}$ as well. Therefore, $\OL\omega_{24}$ is $3$-stiff. Vectors $\pm {\bf e}_i$, $i=\OL {1,4},$ and $\(\pm \frac {1}{2},\pm \frac {1}{2},\pm \frac {1}{2},\pm \frac {1}{2}\)$ form another $24$-cell denoted by~$\W \omega_{24}$.

\begin {table}
\centering
\begin{tabular}{|c|c|c|c|}
\hline
{\bf Configuration $\omega_N$} & {\bf Sphere} & {\bf Assumption} & {\bf Absolute minima}\\
\hline
Regular $2m$-gon & $S^1$ & $f^{(2m-2)}$ is convex & regular $2m$-gon rotated by $\frac{\pi}{2m}$ rad\\
\hline
Regular $(2m+1)$-gon & $S^1$ & $f^{(2m-1)}$ is concave & antipod of regular $(2m+1)$-gon\\
\hline
Regular $d$-simplex $\omega_{d+1}^\ast$ & $S^{d-1}$ & $f'$ is concave & regular $d$-simplex $-\omega^\ast_{d+1}$ \\
\hline
Cross-polytope $\omega^\ast_{2d}$, $d\geq 3$ & $S^{d-1}$ & $f''$ is convex & cube $U_d$ \\
\hline
Cube $U_d$, $d\geq 3$ & $S^{d-1}$ & $f''$ is convex & cross-polytope $\omega^\ast_{2d}$ \\
\hline
$24$-cell $\OL\omega_{24}$ & $S^3$ & $f^{(4)}$ is convex & $24$-cell $\W\omega_{24}$\\
\hline
Schl\"affi $\OL\omega_{27}$ & $S^5$ & $f'''$ is concave & $-\OL\omega_{27}$  \\
\hline
McLaughlin $\OL\omega_{275}$ & $S^{21}$ & $f'''$ is concave & $-\OL\omega_{275}$\\
\hline
\end {tabular}
\caption {Immediate consequences of Theorems \ref {2m-1} and \ref {2m-1w}. We assume that $f:[0,4]\to(-\infty,\infty]$ is continuous on $(0,4]$ with $f(0)=\lim\limits_{t\to 0^+} f(t)$ and differentiable on $(0,4)$ with an additional assumption on $f$ specified on $(0,4)$ in the third column. Then the potential $p_f(\cdot, \omega_N)$ attains its absolute minimum over the specified sphere at points of the set in the fourth column. If the assumption in the third column holds strictly, then no other absolute minima exist.}\label {table4}
\end {table}

\begin {lemma}\label {24-cell}
$\mathcal D_3(\OL\omega_{24})=\W \omega_{24}$.
\end {lemma}
\begin {proof}
The inclusion $\W\omega_{24}\subset \mathcal D_3(\OL \omega_{24})$ is immediate. Let ${\bf x}=(x_1,x_2,x_3,x_4)\in \mathcal D_3(\OL \omega_{24})$. There are at most three distinct numbers among all dot products $\frac {\pm x_i\pm x_j}{\sqrt{2}}$, $1\leq i\neq j\leq 4$. Then all non-zero coordinates of ${\bf x}$ must have the same absolute value (denoted by $a$). If ${\bf x}$ had two or three non-zero coordinates, it would form dot products $0,\pm \frac {a}{\sqrt{2}},\pm \sqrt{2}a$ with vectors from $\OL\omega_{24}$. Therefore, ${\bf x}$ has one non-zero coordinate ($a=1$) or all its coordinates equal $\pm a$ ($a=1/2$); i.e.,~${\bf x}\in \W\omega_{24}$.
\end {proof}

Immediate consequences of Theorems \ref {2m-1} and \ref {2m-1w} are listed in Table \ref {table4}.

\section {Proofs of Theorems \ref {2m-1} and \ref {2m-1w}}\label {proofconst}

We first establish the following interpolation lemma for the potential function.
It is a ``slight" modification of \cite[Lemma 3.1]{BorMinMax}.
The case $\nu=0$ corresponds to stiff configurations and Theorem \ref {2m-1}. The case $\nu=1$ corresponds to strongly sharp configurations and Theorem \ref {2m-1w}.
\begin {lemma}\label {Interp}
Suppose $\nu=0$ or $1$. Let $n\geq 1+\nu$, $L:=2n-2-\nu$, and $-1\leq t_1<\ldots<t_n<1$ be arbitrary nodes such that $t_1>-1$ if $\nu=0$ and $t_1=-1$ if $\nu=1$. Suppose also that $g:[-1,1]\to (-\infty,\infty]$ is a function continuous on $[-1,1)$ with $g(1)=\lim\limits_{t\to 1^-}g(t)$ and differentiable on $(-1,1)$ with the derivative $g^{(L)}$ being convex on $(-1,1)$. Let $q\in \mathbb P_{L+1}$ be the polynomial such that
\begin {equation}\label {i3}
q(t_i)=g(t_i),\ \ i=1,\ldots,n,\ \ \text {and} \ \ q'(t_i)=g'(t_i),\ \ i=1+\nu,\ldots,n.
\end {equation}
Then 
\begin {equation}\label {i5}
q(t)\leq g(t), \ \ t\in [-1,1]. 
\end {equation}
If, in addition, the derivative $g^{(L)}$ is strictly convex on $(-1,1)$, then the inequality in \eqref {i5} is strict for $t\in [-1,1]\setminus \{t_1,\ldots,t_{n}\}$. 
\end {lemma}
\begin {proof}
The assertion of the lemma follows immediately when $L=0$; i.e., when $n=1$ and $\nu=0$. Therefore, we assume that $L\geq 1$ (then $n\geq 2$).
Inequality \eqref {i5} holds trivially for $t=t_1,\ldots,t_{n}$ (as equality) and, if $g(1)=\infty$, at $t=1$ (as a strict inequality). Therefore, we choose arbitrary $x\in [-1,1]\setminus \{t_1,\ldots,t_{n}\}$ assuming additionally that $x<1$ if $g(1)=\infty$. Denote
$$
h(t):=q(t)+b(t-t_1)^{2-\nu}\prod_{i=2}^{n}(t-t_i)^2,
$$
where the constant $b$ is chosen so that $h(x)=g(x)$. Observe that the polynomial $h$ also satisfies interpolation conditions \eqref {i3}. The function $\tau(t):=g(t)-h(t)$ has at least $n+1$ distinct zeros in $[-1,1]$ at least $n-\nu$ of which have multiplicity two and are located in $(-1,1)$. By the Rolle's theorem, the derivative $\tau'$ has at least $2n-\nu=L+2$ distinct zeros in $(-1,1)$, which implies that the derivative $\tau^{(L)}$ has at least three distinct zeros in $(-1,1)$. The polynomial $h^{(L)}$ has degree at most two with the coefficient of $t^2$ being $\frac {(L+2)!}{2}b$. It interpolates the convex function $g^{(L)}$ at at least three distinct points in $(-1,1)$. Then $b\geq 0$. Since $(x-t_1)^{2-\nu}>0$ both for $\nu=0$ and $\nu=1$, we have
\begin {equation}\label {i7}
g(x)=h(x)=q(x)+b(x-t_1)^{2-\nu}\prod_{i=2}^{n}(x-t_i)^2\geq q(x),
\end {equation}
which proves \eqref {i5}.
If $g^{(L)}$ is strictly convex on $(-1,1)$, then $b\neq 0$ (otherwise $g^{(L)}$ would be interpolated by the polynomial $h^{(L)}$ of degree at most one at three distinct points). Therefore, $b>0$ and the inequality in \eqref {i5} is strict on $[-1,1]\setminus \{t_1,\ldots,t_{n}\}$. 
\end {proof}

\begin {lemma}\label {w2m-1}
Let $\omega_N$ be an $m$-stiff configuration on $S^d$, $d,m\geq 1$. Then there is no point in $\mathcal D_m(\omega_N)$ that forms dot products $1$ or $-1$ with points from $\omega_N$. Let $\omega_N$ be strongly $m$-sharp. Then $\omega_N$ contains no antipodal pair and $\mathcal D_m(\omega_N)=\emptyset$.
\end {lemma}
\begin {proof}
Assume to the contrary that there is a point ${\bf y}\in \mathcal D_m(\omega_N)$ such that ${\bf y}\cdot {\bf x}_\ell=1$ or $-1$ for some ${\bf x}_\ell\in \omega_N$. Let $-1\leq t_1<\ldots<t_k\leq 1$, $k\leq m$, be the distinct dot products ${\bf y}$ forms with vectors from $\omega_N$. If ${\bf y}\cdot {\bf x}_\ell=1$, then $t_k=1$ and we let 
$p(t):=(t-t_1)^2\cdots(t-t_{k-1})^2(1-t)$. If ${\bf y}\cdot {\bf x}_\ell=-1$, then $t_1=-1$ and we let $p(t):=(t+1)(t-t_2)^2\cdots(t-t_k)^2$. Define $u({\bf x}):=p({\bf x}\cdot {\bf y})$, ${\bf x}\in \RR^{d+1}$. Then $u$ has degree at most $2m-1$, $u({\bf x}_i)=0$, $i=1,\ldots,N$, and $u({\bf x})>0$ for $\sigma_d$-almost all ${\bf x}\in S^d$. Then, since $\omega_N$ is a $(2m-1)$-design,
$$
0=\frac {1}{N}\sum\limits_{i=1}^{N}u({\bf x}_i)=\int_{S^d} u({\bf x})\ \! d\sigma_d({\bf x})>0.
$$
This contradiction proves the first part of the lemma.

Let now $\omega_N$ be strongly $m$-sharp. Assume to the contrary that $\omega_N$ contains an antipodal pair ${\bf x}_\ell=-{\bf x}_i$. Let $-1=t_1<t_2<\ldots<t_m<1$ be dot products between distinct points in $\omega_N$. Letting $p(t):=(t+1)(t-t_2)^2\cdots(t-t_m)^2(1-t)$ and $u({\bf x}):=p({\bf x}\cdot {\bf x}_\ell)$ and using a similar argument, we arrive at a contradiction proving that $\omega_N$ has no antipodal pair.

Finally, assume to the contrary that there is a point ${\bf z}\in \mathcal D_m(\omega_N)$. Let $-1\leq t_1<\ldots<t_k\leq 1$, $k\leq m$, be distinct elements in $D({\bf z},\omega_N)$. Letting $p(t):=(t-t_1)^2\cdots(t-t_k)^2$ and $u({\bf x}):=p({\bf x}\cdot {\bf z})$ and using a similar argument again, we arrive at a contradiction proving that $\mathcal D_m(\omega_N)$ is empty.
\end {proof}

Replacing $g$ and $q$ with $-g$ and $-q$ in \cite[Lemma 3.5]{BorMinMax} we obtain the following.

\begin {lemma}\label {gen}
Suppose that $d,n\geq 1$ and that a function $g:[-1,1]\to (-\infty,\infty]$ and a set $A$ consisting of numbers $-1\leq t_1<\ldots<t_{n}\leq 1$ are arbitrary. Suppose also that $q$ is a polynomial such that $q(t_i)=g(t_i)$, $i=1,\ldots,n$, and
\begin {equation}\label {1geq}
q(t)\leq g(t),\ \ t\in [-1,1].
\end {equation}
If $\omega_N=\{{\bf x}_1,\ldots,{\bf x}_N\}\subset S^{d}$ is a spherical $M$-design, where $M\geq {\rm deg}\ \! q$, and ${\bf y}\in S^{d}$ is any point such that $D({\bf y},\omega_N)\subset A$, then the potential 
$$
p^g({\bf x},\omega_N):=\sum\limits_{i=1}^{N}g({\bf x}\cdot {\bf x}_i), \ \ \ {\bf x}\in S^{d},
$$
attains its absolute minimum over $S^{d}$ at point ${\bf y}$. 

If, in addition, the inequality in \eqref {1geq} is strict for every $t\in [-1,1]\setminus A$, then any point ${\bf z}\in S^{d}$ such that $D({\bf z},\omega_N)\not\subset A$ is not a point of absolute minimum of $p^g(\ \!\cdot\ \!,\omega_N)$ on~$S^{d}$.
\end {lemma}

We are now ready to prove Theorems \ref {2m-1} and \ref {2m-1w}.

\begin {proof}[Proof of Theorem \ref {2m-1}]
Let ${\bf y}$ be any point from $\mathcal D_m(\omega_N)$. Let $-1<\tau_1<\ldots<\tau_m<1$ be distinct elements of the set $D({\bf y},\omega_N)=\{{\bf y}\cdot{\bf x}_i : 1\leq i\leq N\}$. There are exactly $m$ dot products there in view of Remark \ref {m} with $\tau_1>-1$ and $\tau_m<1$ in view of Lemma~\ref {w2m-1}. We apply Lemma~\ref {Interp} with $\nu=0$, $n=m$, and $g(t)=f(2-2t)$ so that $g^{(L)}$ is convex on $(-1,1)$, where $L=2m-2$. Then the Hermite interpolating polynomial $q\in \mathbb P_{2m-1}$ for $g$ at the nodes $\tau_i$ satisfies $q(t)\leq g(t)$, $t\in [-1,1]$. Since $\omega_N$ is a spherical ($2m-1$)-design and $D({\bf y},\omega_N)= \{\tau_1,\ldots,\tau_m\}$, by Lemma \ref {gen} (with $n=m$ and $M=2m-1$), the potential $p_f(\cdot,\omega_N)=p^g(\cdot,\omega_N)$ attains its absolute minimum over $S^{d}$ at~${\bf y}$. 

If $f^{(2m-2)}$ is strictly convex on $(0,4)$, then so is $g^{(2m-2)}$ on $(-1,1)$ and, by Lemma \ref {Interp}, $q(t)<g(t)$, $t\in [-1,1]\setminus \{\tau_1,\ldots,\tau_m\}$. If ${\bf z}\in S^{d}$ is any point not in $\mathcal D_m(\omega_N)$, then the set $D({\bf z},\omega_N)$ contains at least $m+1$ distinct dot products. Hence, $D({\bf z},\omega_N)\not\subset \{\tau_1,\ldots,\tau_m\}$. In view of Lemma \ref {gen}, the point ${\bf z}$ is not a point of absolute minimum of $p^g(\cdot,\omega_N)$ on $S^{d}$, and, hence, of $p_f(\cdot,\omega_N)$. Thus, $\mathcal D_m(\omega_N)$ is the set of all points of absolute minimum of $p_f(\cdot,\omega_N)$ on $S^{d}$.  
\end {proof}

\begin {proof}[Proof of Theorem \ref {2m-1w}] Let $-1< \tau_1<\ldots<\tau_m<\tau_{m+1}=1$ be the distinct values of the dot product between points in $\omega_N$. Here, we have $\tau_1>-1$ in view of Lemma \ref {w2m-1}. We apply Lemma \ref {Interp} with $\nu=1$, $n=m+1$, $t_1=-\tau_{m+1}=-1$, $t_2=-\tau_m$, ..., $t_{m+1}=-\tau_1$, and $g(t)=f(2-2t)$ so that $g^{(L)}$ is convex on $(-1,1)$, where $L=2m-1$. Then the polynomial $q\in \mathbb P_{2m}$, which interpolates both $g$ and $g'$ at $t_2,\ldots,t_{m+1}$ and only $g$ at $t_1$ satisfies $q(t)\leq g(t)$, $t\in [-1,1]$. Choose any point ${\bf y}$ in $-\omega_N$. Then $D({\bf y},\omega_N)\subset \{t_1,\ldots,t_{m+1}\}$. Since $\omega_N$ is a $2m$-design, by Lemma \ref {gen} (with $n=m+1$ and $M=2m$) the potential $p_f(\cdot,\omega_N)=p^g(\cdot,\omega_N)$ attains its absolute minimum over $S^d$ at ${\bf y}$.

Let now $f^{(2m-1)}$ be strictly concave on $(0,4)$. Then $g^{(2m-1)}$ is strictly convex on $(-1,1)$. Let ${\bf z}$ be an arbitrary point in $S^d\setminus (-\omega_N)$. Assume to the contrary that $D({\bf z},\omega_N)\subset \{t_1,\ldots,t_{m+1}\}$. If $D({\bf z},\omega_N)$ did not contain $t_1=-1$, then it would contain at most $m$ distinct elements making the set $\mathcal D_m(\omega_N)$ non-empty. By Lemma \ref {w2m-1}, this is not possible. Consequently, $D({\bf z},\omega_N)$ contains $-1$. Then ${\bf z}\in -\omega_N$ contradicting its choice. Thus, $D({\bf z},\omega_N)\not\subset \{t_1,\ldots,t_{m+1}\}$. By Lemma~\ref {Interp}, we have $q(t)<g(t)$, $t\in [-1,1]\setminus \{t_1,\ldots,t_{m+1}\}$. Then by Lemma \ref {gen}, the point ${\bf z}$ is not a point of absolute minimum of the potential $p_f(\cdot,\omega_N)=p^g(\cdot,\omega_N)$ on~$S^d$.
\end {proof}

\section {Dot product distribution related to stiff configurations}\label {DDP}

In this section, we study the distribution of dot products between vectors from $\mathcal D_m(\omega_N)$ and $\omega_N$. 
Spherical designs can be characterized using Gegenbauer polynomials as follows.
\begin {theorem}\label {cd}(see \cite {DelGoeSei1977} or \cite[Theorems 5.2.2 and 5.4.2]{BorHarSafbook})
Let $d,n\geq 1$ and $\omega_N=\{{\bf x}_1,\ldots,{\bf x}_N\}$ be a point configuration on $S^d$. The following are equivalent:
\begin {itemize}
\item[(i)]
$\omega_N$ is a spherical $n$-design;
\item [(ii)]
$
\sum\limits_{i=1}^N\sum\limits_{j=1}^{N}P_k^{(d)}({\bf x}_i\cdot{\bf x}_j)=0,\ \ k=1,\ldots,n;
$
\item[(iii)]
For every polynomial $q\in \mathbb P_n$, we have $p^q({\bf y},\omega_N):=\sum_{i=1}^{N}q({\bf y}\cdot{\bf x}_i)=a_0(q)N$, ${\bf y}\in S^d$, where
$$
a_0(q):=\int_{-1}^{1}q(t)w_d(t)\ \!dt
$$
is the $0$-th Gegenbauer coefficient of polynomial $q$.
\end {itemize}
\end {theorem}
One can use Theorem \ref {cd} part (ii) to find the strength of every spherical design mentioned in this paper. 
We next restate Remark \ref {m1} formally and, for completeness, provide its proof.
Here, $\varphi_1,\ldots,\varphi_m$ are the fundamental polynomials for set of nodes \eqref {fund}; that is, $\varphi_i\in \mathbb P_{m-1}$, $\varphi_i(\kappa_i^m)=1$, and $\varphi_i(\kappa_j^m)=0$, $j\neq i$, $i=1,\ldots,m$. 

\begin {proposition}\label {mtight}
If $\omega_N$ is an $m$-stiff configuration on $S^d$, then for every ${\bf z}\in \mathcal D_m(\omega_N)$, the set $D({\bf z},\omega_N)$ contains exactly $m$ distinct elements located in $(-1,1)$, which are the zeros $\kappa_1^{m},\ldots,\kappa_{m}^{m}$ of Gegenbauer polynomial $P_m^{(d)}$. Furthermore, the number of indices $i$ such that ${\bf z}\cdot {\bf x}_i=\kappa_j^{m}$ does not depend on ${\bf z}$ and equals $a_0(\varphi_j) N$, $j=1,\ldots,m$.

In particular, if $m=2$, then for every ${\bf z}\in \mathcal D_2(\omega_N)$, we have $D({\bf z},\omega_N)=\left\{-\frac {1}{\sqrt{d+1}},\frac {1}{\sqrt{d+1}}\right\}$, and, if $m=3$, then for every ${\bf z}\in \mathcal D_3(\omega_N)$, we have $D({\bf z},\omega_N)=\left\{\pm \sqrt{\frac {3}{d+3}},0\right\}$.
\end {proposition}
\begin {proof}
Denote by $-1\leq t_1<\ldots<t_m\leq 1$ the distinct elements in $D({\bf z},\omega_N)$ (in view of Remark \ref {m}, there are exactly $m$ of them). Let $M_j:=\#\{i : {\bf z}\cdot {\bf x}_i=t_j\}$, $j=1,\ldots,m$. Denote $P(t):=(t-t_1)\cdots (t-t_m)$. Since $\omega_N$ is a $(2m-1)$-design, using item (iii) of Theorem \ref {cd}, for every polynomial $p\in \mathbb P_{m-1}$, we have
$$
a_0(pP)=\frac {1}{N}\sum\limits_{i=1}^N p({\bf z}\cdot {\bf x}_i)P({\bf z}\cdot {\bf x}_i)=\frac {1}{N}\sum\limits_{j=1}^{m}M_jp(t_j)P(t_j)=0.
$$
Then $P\bot \ \!\mathbb P_{m-1}$ on $[-1,1]$ with weight $w_d$; that is, $P$ is a constant multiple of the Gegenbauer polynomial $P_m^{(d)}$. Then $t_j=\kappa_j^m$, $j=1,\ldots,m$, and we have
$$
a_0(\varphi_k)=\frac {1}{N}\sum\limits_{i=1}^{N}\varphi_k({\bf z}\cdot{\bf x}_i)=\frac {1}{N}\sum\limits_{j=1}^{m}M_j\varphi_k(\kappa_j^m)=\frac {M_k}{N}; 
$$
that is, $M_k=a_0(\varphi_k)N$, $k=1,\ldots,m$.

Since $P_2^{(d)}(t)=\frac {d+1}{d}t^2-\frac {1}{d}$ and $P_3^{(d)}(t)=\frac {d+3}{d}t^3-\frac {3}{d}t$, the sets of zeros of these polynomials give the set $D({\bf z},\omega_N)$ for $m=2$ and $3$.
\end {proof}

We conclude this section with the following consequence of Lemma \ref {w2m-1}.
\begin {proposition}\label {sst}
The classes of stiff and strongly sharp configurations on $S^{d}$, $d\geq 1$, are disjoint. \end {proposition}

\begin {proof}[Proof of Proposition \ref {sst}]
Assume to the contrary that there exists a configuration $\omega_N\subset S^d$ which is $k$-stiff and strongly $m$-sharp for some $k,m\geq 1$. By Lemma~\ref {w2m-1}, the set $\mathcal D_m(\omega_N)$ is empty. Then $k\geq m+1$. Since $\omega_N$ is $k$-stiff, it must be at least a $(2m+1)$-design. If $\omega_N$ is not $(m+1)$-stiff, then, by definition, $\mathcal D_{m+1}(\omega_N)=\emptyset$. If $\omega_N$ is $(m+1)$-stiff, then by Lemma \ref {w2m-1}, the set $\mathcal D_{m+1}(\omega_N)$ does not contain any point from $\omega_N$. However, since $\omega_N$ is strongly $m$-sharp, every point from $\omega_N$ is in $\mathcal D_{m+1}(\omega_N)$. This contradiction proves the proposition.
\end {proof}

\section {Symmetrized regular simplex and its dual}\label {symm}

In this section, we describe universal minima for the first pair of stiff configurations.
Let $d\geq 3$ be odd. Denote $\OL\omega_{2d+2}=\omega^\ast_{d+1}\cup(-\omega^\ast_{d+1})$, where $\omega^\ast_{d+1}=\{{\bf x}_0,{\bf x}_1\ldots,{\bf x}_{d}\}$ is the set of vertices of a regular $d$-simplex inscribed in $S^{d-1}$. Since $\omega_{d+1}^\ast$ is a spherical $2$-design, the symmetrized regular simplex $\OL\omega_{2d+2}$ is a spherical $3$-design. Denote $\kappa(d):=\binom {d+1}{\frac {d+1}{2}}$ and let
\begin{equation}\label{SymSimp}
\OL\omega_{\kappa(d)}:=\left\{\frac {2\sqrt{d}}{d+1}\sum\limits_{i\in I}{\bf x}_i : I\subset \{0,1,\ldots,d\},\ \# I=\frac {d+1}{2}\right\}.
\end{equation}
It is not difficult to see that $\OL\omega_{\kappa(d)}\subset S^{d-1}$. By Lemma \ref {kappa(d)} below, $\mathcal D_2(\OL\omega_{2d+2})=\OL\omega_{\kappa(d)}$. Thus, $\OL\omega_{2d+2}$ is $2$-stiff.  The configuration $\OL\omega_{\kappa(d)}$, in turn, is a spherical $3$-design (see Lemma~\ref {omegakappa} below) with $\mathcal D_2(\OL\omega_{\kappa(d)})=\OL\omega_{2d+2}$ in view of Lemma \ref {omegakappa1}. Then $\OL\omega_{\kappa(d)}$ is also $2$-stiff. We remark that for $d=2$, the configuration $\OL\omega_{2d+2}$ is a regular hexagon on $S^1$, which is $3$-stiff, and its dual is also a regular hexagon. 

\begin {proposition}\label {ns}
For $d\geq 4$ even, $\mathcal D_2(\OL\omega_{2d+2})=\emptyset$ and $\OL\omega_{2d+2}$ is not stiff. 
\end {proposition}
\begin {proof}
The configuration $\OL\omega_{2d+2}$ is a $3$-design. Assume to the contrary that $\mathcal D_2(\OL\omega_{2d+2})\neq \emptyset$. Then $\OL\omega_{2d+2}$ is $2$-stiff. By Proposition \ref {mtight}, there is a point ${\bf x}\in S^{d-1}$ such that ${\bf x}\cdot {\bf x}_i=\pm \frac {1}{\sqrt{d}}$, $i=0,1,\ldots,d$. Let $k$ be the number of indices $i$ such that ${\bf x}\cdot {\bf x}_i=\frac {1}{\sqrt{d}}$. Then, since $d+1$ is odd,
$$
{\bf x}\cdot \sum_{i=0}^{d}{\bf x}_i=\sum\limits_{i=0}^{d}{\bf x}\cdot {\bf x}_i=\frac {k}{\sqrt{d}}-\frac {d+1-k}{\sqrt{d}}=\frac {2k-(d+1)}{\sqrt{d}}\neq 0
$$
contradicting the fact that the center of mass of $\omega_{d+1}^\ast$ is at the origin. Therefore, $\mathcal D_2(\OL\omega_{2d+2})=\emptyset$, and $\OL\omega_{2d+2}$ cannot be $1$- or $2$-stiff. Since $d\geq 4$, $\OL\omega_{2d+2}$ is not a $4$-design. Then $\OL\omega_{2d+2}$ cannot be $m$-stiff for any $m\geq 3$.
\end {proof}

\begin {lemma}\label {kappa(d)}
For every $d\geq 3$ odd, 
$\mathcal D_2(\OL\omega_{2d+2})=\OL\omega_{\kappa(d)}$.
\end {lemma}
\begin {proof}
Let ${\bf x}_0,{\bf x}_1,\ldots,{\bf x}_d$ be vertices of a regular $d$-simplex on $S^{d-1}$ such that $\OL\omega_{2d+2}=\{\pm {\bf x}_i : i=0,1,\ldots,d\}$.
The inclusion $\OL\omega_{\kappa(d)}\subset \mathcal D_2(\OL\omega_{2d+2})$ holds. Indeed, choose any $\frac{d+1}{2}$-element subset $I\subset \{0,1,\ldots,d\}$. If $k$ is an index in $I$, then
\begin {equation}\label {pmd}
\(\frac {2\sqrt{d}}{d+1}\sum\limits_{i\in I}{\bf x}_i\)\cdot {\bf x}_k= \frac {2\sqrt{d}}{d+1}\(1-\frac {1}{d}\cdot \frac{d-1}{2}\)=\frac {1}{\sqrt{d}}.
\end {equation}
If $k$ is an index not in $I$, then
\begin {equation}\label {pnd}
\(\frac {2\sqrt{d}}{d+1}\sum\limits_{i\in I}{\bf x}_i\)\cdot {\bf x}_k=\frac {2\sqrt{d}}{d+1}\cdot \(-\frac {1}{d}\cdot \frac {d+1}{2}\)=- \frac {1}{\sqrt{d}}.
\end {equation}
Then every point from $\OL\omega_{\kappa(d)}$ forms only two distinct values of the dot product with points from $\OL\omega_{2d+2}$ (which are $\pm \frac {1}{\sqrt{d}}$) and, hence, is contained in $\mathcal D_2(\OL\omega_{2d+2})$.

To show the opposite inclusion, choose any ${\bf y}\in \mathcal D_2(\OL\omega_{2d+2})$. Since points ${\bf x}_0,{\bf x}_1,\ldots,{\bf x}_d$ form a $2$-design, the set $\OL\omega_{2d+2}$ is a $3$-design. Since $\mathcal D_2(\OL\omega_{2d+2})\neq \emptyset$ (it contains $\OL\omega_{\kappa(d)}$), the set $\OL\omega_{2d+2}$ is $2$-stiff. By Proposition \ref {mtight}, we then have ${\bf y}\cdot{\bf x}_i=\pm \frac {1}{\sqrt{d}}$, $i=0,1,\ldots,d$. Since $d+1$ is even and 
$$
\sum_{i=0}^{d}{\bf y}\cdot{\bf x}_i={\bf y}\cdot\sum_{i=0}^{d}{\bf x}_i={\bf y}\cdot{\bf 0}=0,
$$ 
exactly half of dot products ${\bf y}\cdot{\bf x}_i$ equal $\frac {1}{\sqrt{d}}$ and the other half equals $-\frac {1}{\sqrt{d}}$. Let $I\subset \{0,1,\ldots,d\}$ be the $\frac {d+1}{2}$-element subset such that ${\bf y}\cdot{\bf x}_i=\frac {1}{\sqrt{d}}$, $i\in I$. Pick some index $m\in I$ and some index $k\in J:=\{0,1,\ldots,d\}\setminus I$. Then for every $i\in I\setminus \{m\}$, we have ${\bf y}\cdot({\bf x}_i-{\bf x}_m)=0$ and, for every $i\in J\setminus \{k\}$, we have ${\bf y}\cdot({\bf x}_i-{\bf x}_k)=0$. Let us show that the set $$
T:=\{{\bf x}_i-{\bf x}_m\}_{i\in I\setminus\{m\}}\cup \{{\bf x}_i-{\bf x}_k\}_{i\in J\setminus\{k\}}
$$
of $d-1$ vectors in $\RR^d$ (orthogonal to ${\bf y}$) is linearly independent. Let $c_i$, $i\in \{0,1,\ldots,d\}\setminus\{m,k\}$ be any numbers such that
$$
\sum\limits_{i\in I\setminus \{m\}}c_i({\bf x}_i-{\bf x}_m)+\sum\limits_{i\in J\setminus \{k\}}c_i({\bf x}_i-{\bf x}_k)={\bf 0}.
$$
Then 
\begin {equation}\label {cij}
\sum\limits_{i=0\atop i\neq m,k}^{d}c_i({\bf x}_i-{\bf x}_m)-\(\sum\limits_{i\in J\setminus \{k\}}c_i\)({\bf x}_k-{\bf x}_m)={\bf 0}.
\end {equation}
Since ${\bf x}_0,\ldots,{\bf x}_d$ form a non-degenerate simplex in $\RR^d$, the set of vectors ${\bf x}_i-{\bf x}_m$, where $i\in \{0,1,\ldots,d\}\setminus \{m\}$ is linearly independent. Then \eqref {cij} implies that $c_i=0$, $i\neq m,k$; that is, $T$ is linearly independent. Thus, ${\bf y}$ appears to be orthogonal to the $(d-1)$-dimensional subspace spanned by $T$. Then ${\bf y}$ belongs to the intersection of the one-dimensional subspace $T^\bot$ with $S^{d-1}$; that is there exist only two (opposite) vectors that ${\bf y}$ can equal to. 

Let ${\bf u}:=\frac {2\sqrt{d}}{d+1}\sum_{j\in I}{\bf x}_j$. Then ${\bf u}\in S^{d-1}\cap \OL\omega_{\kappa(d)}$ and equation \eqref {pmd} implies that, for every $i\in I\setminus \{m\}$, we have ${\bf u}\cdot({\bf x}_i-{\bf x}_m)=\frac {1}{\sqrt{d}}-\frac {1}{\sqrt{d}}=0$. Equation \eqref {pnd} implies that, for every $i\in J\setminus \{k\}$, we have ${\bf u}\cdot ({\bf x}_i-{\bf x}_k)=-\frac {1}{\sqrt{d}}+\frac {1}{\sqrt{d}}=0$. Thus, ${\bf u}\in T^\bot \cap S^{d-1}$. This implies that ${\bf y}={\bf u}$ or ${\bf y}=-{\bf u}$. Since $\sum_{i=0}^{d}{\bf x}_i={\bf 0}$, we have $-{\bf u}=\frac {2\sqrt{d}}{d+1}\sum_{j\in J}{\bf x}_j\in \OL\omega_{\kappa(d)}$. Thus, ${\bf y}\in \OL\omega_{\kappa(d)}$.
\end {proof}

\begin {lemma}\label {omegakappa}
For every $d\geq 3$ odd, the configuration $\OL\omega_{\kappa(d)}\subset S^{d-1}$ is a spherical $3$-design. 
\end {lemma}
\begin {proof}
To prove the lemma, we will use Theorem \ref {cd}.
Since $\OL\omega_{\kappa(d)}$ is centrally symmetric and Gegenbauer polynomials $P_n^{(d-1)}$, $n=1,3$, are odd, we have
\begin {equation}\label {symmodd}
M_n^{(d-1)}(\OL\omega_{\kappa(d)}):=\sum\limits_{{\bf x}\in \OL\omega_{\kappa(d)}}\sum\limits_{{\bf y}\in \OL\omega_{\kappa(d)}}P_n^{(d-1)}({\bf x}\cdot{\bf y})=0
\end {equation}
for $n=1,3$. It remains to show equality \eqref {symmodd} for $n=2$.
Let ${\bf x}_0,{\bf x}_1,\ldots,{\bf x}_d$ be the vertices of the regular simplex in the definition of the configuration $\OL\omega_{\kappa(d)}$ in \eqref {SymSimp}. Let ${\bf x}_I:=\frac {2\sqrt{d}}{d+1}\sum_{i\in I}{\bf x}_i$, where $I$ is any $\ell:=\frac {d+1}{2}$-element subset of $\{0,1,\ldots,d\}$. Let $0\leq k\leq \ell$ be any integer and $J\subset \{0,1,\ldots,d\}$ be any $\ell$-element subset whose intersection with $I$ has exactly $k$ elements. There are $\binom {\ell}{k}\binom {\ell}{\ell-k}=\binom{\ell}{k}^2$ such sets $J$ and, for every such $J$, 
\begin {equation*}
\begin {split}
{\bf x}_I\cdot {\bf x}_J&=\frac {4d}{(d+1)^2}\(\sum_{i\in I\cap J}{\bf x}_i+\sum_{i\in I\setminus J}{\bf x}_i\)\cdot \sum_{j\in J}{\bf x}_j\\
&=\frac {d}{\ell^2}\(\sum_{i\in I\cap J}\(1-\frac {\ell-1}{d}\)+\sum_{i\in I\setminus J}\(-\frac {\ell}{d}\)\)=\frac {2k}{\ell}-1,
\end {split}
\end {equation*}
$k=0,1,\ldots,\ell$. Then
\begin {equation}\label {M2}
\begin {split}
M_2^{(d-1)}&(\OL\omega_{\kappa(d)})=\kappa(d)\sum\limits_{k=0}^{\ell}\binom{\ell}{k}^2P_2^{(d-1)}\(\frac {2k}{\ell}-1\)\\
&=\frac {d\kappa(d)}{d-1}\sum\limits_{k=0}^{\ell}\binom{\ell}{k}^2\(\(\frac {2k}{\ell}-1\)^2-\frac {1}{d}\)\\
&=\frac {d\kappa(d)}{d-1}\(\frac {4}{\ell^2}\sum\limits_{k=0}^{\ell}k^2\binom{\ell}{k}^2-\frac {4}{\ell}\sum\limits_{k=0}^{\ell}k\binom{\ell}{k}^2+\frac {d-1}{d}\sum\limits_{k=0}^{\ell}\binom{\ell}{k}^2\).
\end {split}
\end {equation}
Using the standard formula $\sum\limits_{k=0}^{\ell}\binom{\ell}{k}^2=\binom{2\ell}{\ell}$ (which can be seen by counting the total number of subsets $J$ above), we obtain that
$$
\sum\limits_{k=0}^{\ell}k^2\binom{\ell}{k}^2=\sum\limits_{k=1}^{\ell}\ell^2\binom{\ell-1}{k-1}^2=\ell^2\sum\limits_{k=0}^{\ell-1}\binom{\ell-1}{k}^2=\ell^2\binom{2\ell-2}{\ell-1}.
$$
Finally, from \eqref {symmodd} we obtain that
\begin {equation*}
\begin{split}
&0=\frac {M_1^{(d-1)}(\OL\omega_{\kappa(d)})}{\kappa(d)}=\sum\limits_{k=0}^{\ell}\binom{\ell}{k}^2P_1^{(d-1)}\(\frac {2k}{\ell}-1\)\\
&=\sum\limits_{k=0}^{\ell}\binom{\ell}{k}^2\(\frac {2k}{\ell}-1\)=\frac {2}{\ell}\sum\limits_{k=0}^{\ell}k\binom{\ell}{k}^2-\binom{2\ell}{\ell}.
\end {split}
\end {equation*}
Then $\sum_{k=0}^{\ell}k\binom{\ell}{k}^2=\frac {\ell}{2}\binom{2\ell}{\ell}$ and from \eqref {M2} we obtain
$$
M_2^{(d-1)}(\OL\omega_{\kappa(d)})=\frac {d\kappa(d)}{d-1}\(4\binom{2\ell-2}{\ell-1}-2\binom{2\ell}{\ell}+\frac {2\ell-2}{2\ell-1}\binom{2\ell}{\ell}\)=0.
$$
Thus, $\OL\omega_{\kappa(d)}$ is a spherical $3$-design.
\end {proof}

\begin {lemma}\label {omegakappa1}
For every $d\geq 3$ odd, $\mathcal D_2(\OL\omega_{\kappa(d)})=\OL\omega_{2d+2}$.
\end {lemma}
\begin {proof}
Relations \eqref {pmd} and \eqref {pnd} imply that $\OL\omega_{2d+2}\subset \mathcal D_2(\OL\omega_{\kappa(d)})$.
Let ${\bf x}$ be any element of $\mathcal D_2(\OL\omega_{\kappa(d)})$ and ${\bf x}_0,{\bf x}_1,\ldots,{\bf x}_d$ be the vertices of the regular simplex in the definition of the configuration $\OL\omega_{\kappa(d)}$ in \eqref {SymSimp}. There exist numbers $\alpha_0,\alpha_1,\ldots,\alpha_d$ such that $\sum_{i=0}^{d}\alpha_i=1$ and ${\bf x}=\sum_{i=0}^{d}\alpha_i{\bf x}_i$.
Let $I$ be any $(d+1)/2$-element subset of the set $\{0,1,\ldots,d\}$. Then relations \eqref {pmd} and \eqref {pnd} imply that 
\begin{equation*}
\begin {split}
\theta:=\(\frac {2\sqrt{d}}{d+1}\sum\limits_{i\in I}{\bf x}_i\)\cdot {\bf x}&=\(\frac {2\sqrt{d}}{d+1}\sum\limits_{i\in I}{\bf x}_i\)\cdot\(\sum\limits_{j\in I}\alpha_j{\bf x}_j+\sum\limits_{j\notin I}\alpha_j{\bf x}_j\)
\\&=\frac {1}{\sqrt{d}}\(\sum\limits_{j\in I}\alpha_j-\sum\limits_{j\notin I}\alpha_j\).
\end {split}
\end {equation*}
The set $\OL\omega_{\kappa(d)}$ is a $3$-design by Lemma \ref {omegakappa}. Since $\mathcal D_2(\OL\omega_{\kappa(d)})\neq \emptyset$ (it contains $\OL\omega_{2d+2}$), the set $\OL\omega_{\kappa(d)}$ is $2$-stiff. Then in view of Proposition \ref {mtight}, the above dot product is $\theta=\pm \frac {1}{\sqrt{d}}$. Then $\sum_{i\in I}\alpha_i-\sum_{i\notin I}\alpha_i=\pm 1$. Since $\sum_{i=0}^{d}\alpha_i=1$, for any $I$, we have $\sum_{i\in I}\alpha_i=1$ or $\sum_{i\in I}\alpha_i=0$. 

Since ${\bf x}$ and each ${\bf x}_i$ lie on the unit sphere, we have
\begin {equation*}
\begin {split}
\left|{\bf x}\right|^2&=\(\sum\limits_{i=0}^{d}\alpha_i{\bf x}_i\)^2=\sum\limits_{i=0}^{d}\alpha_i^2-\frac {1}{d}\sum\limits_{i\neq j}\alpha_i\alpha_j=\(1+\frac {1}{d}\)\sum\limits_{i=0}^{d}\alpha_i^2-\frac {1}{d}\(\sum\limits_{i=0}^{d}\alpha_i\)^2\\
&=\(1+\frac {1}{d}\)\sum\limits_{i=0}^{d}\alpha_i^2-\frac {1}{d}=1,
\end {split}
\end {equation*}
which implies that $\sum_{i=0}^{d}\alpha_i^2=1$.

Choose arbitrary indices $j\neq k$ and let $J\subset \{0,1,\ldots,d\}$ be any $(d+1)/2$-element subset such that $k\in J$ and $j\notin J$. Let $J_1:=(J\cup \{j\})\setminus \{k\}$. Since each sum $\sum_{i\in J}\alpha_i$ and $\sum_{i\in J_1}\alpha_i$ equals $0$ or $1$, the difference
$$
\alpha_j-\alpha_k=\sum_{i\in J_1}\alpha_i-\sum_{i\in J}\alpha_i
$$
can only equal $1$, $-1$, or $0$. Let $\alpha:=\max\limits_{i=\OL{0,d}}\alpha_i$.  If it were that $\alpha_i=\alpha$ for all indices $i$, then $\alpha_i$ would equal $\frac {1}{d+1}$ for all $i$ contradicting the equality of $\sum_{i\in J}\alpha_i$ to $0$ or $1$. Thus, for some index $\ell$, we have $\alpha_\ell<\alpha$. For every index $i$, $\alpha-\alpha_i$ can only equal $0$ or $1$ (since $\alpha$ is some maximal $\alpha_j$). Then $\alpha_i=\alpha$ or $\alpha_i=\alpha-1$ for every index $i$. Let $n$ be the number of indices $i$ such that $\alpha_i=\alpha$. Since the sums of $\alpha_i$'s and of their squares both equal $1$, we obtain the system of equations
$$
\begin {cases}
n\alpha+(d+1-n)(\alpha-1)=1, \cr
n\alpha^2+(d+1-n)(\alpha-1)^2=1. \cr
\end {cases}
$$

Multiplying the first equation by $\alpha$ and subtracting it from the second equation, we obtain that $(d+1-n)(1-\alpha)=1-\alpha$. This yields two solutions: $\alpha=1$, $n=1$ and $\alpha=\frac {2}{d+1}$, $n=d$. In the first case, we have $\alpha_\ell=1$ for some index $\ell$ and $\alpha_i=0$ for $i\neq \ell$. Then ${\bf x}={\bf x}_\ell\in \OL\omega_{2d+2}$. In the second case, we have $\alpha_\ell=\frac {1-d}{d+1}$ for some index $\ell$ and $\alpha_i=\frac {2}{d+1}$, $i\neq \ell$. Then ${\bf x}=\frac {2}{d+1}\sum\limits_{i: i\neq \ell}{\bf x}_i+\frac {1-d}{d+1}{\bf x_\ell}$. Using the equality $\sum_{i=0}^{d}{\bf x}_i={\bf 0}$, we obtain that ${\bf x}=-{\bf x}_\ell\in \OL\omega_{2d+2}$. Thus, $\mathcal D_2(\OL\omega_{\kappa(d)})\subset \OL\omega_{2d+2}$ and the assertion of the lemma follows.
\end {proof}

We are ready to prove the main result of this section.

\begin {theorem}\label {symmetrizedsimplex}
Let $d\geq 3$ be odd and $f:[0,4]\to (-\infty,\infty]$ be a function continuous on $(0,4]$ with $f(0)=\lim\limits_{t\to 0^+}f(t)$ and a convex derivative $f''$ on $(0,4)$. Then
\begin {itemize}
\item[(i)]
the potential $p_f(\cdot,\OL\omega_{2d+2})$ of the symmetrized regular $d$-simplex $\OL\omega_{2d+2}$
attains its absolute minimum over $S^{d-1}$ at every point of the set $\OL\omega_{\kappa(d)}$;

\item [(ii)]
the potential $p_f(\cdot,\OL\omega_{\kappa(d)})$ of the configuration $\OL\omega_{\kappa(d)}$
attains its absolute minimum over $S^{d-1}$ at every point of the set $\OL\omega_{2d+2}$.
\end {itemize}
If the convexity of $f''$ is strict on $(0,4)$, then the potential $p_f(\cdot,\OL\omega_{2d+2})$ has no other points of absolute minimum on $S^{d-1}$ in (i) and the potential $p_f(\cdot,\OL\omega_{\kappa(d)})$ has no other points of absolute minimum on $S^{d-1}$ in (ii).
\end {theorem}

\begin {proof}[Proof of Theorem \ref {symmetrizedsimplex}]
The set $\OL\omega_{2d+2}$ is a $3$-design as a symmetrization about the origin of a $2$-design. 
In view of Lemma \ref {kappa(d)}, we have $\mathcal D_2(\OL\omega_{2d+2})=\OL\omega_{\kappa(d)}\neq \emptyset$; that is, $\OL\omega_{2d+2}$ is $2$-stiff. Theorem~\ref {2m-1} now implies that the potential $p_f(\ \! \cdot\ \!,\OL\omega_{2d+2})$ attains its absolute minimum on $S^{d-1}$ at every point of $\OL\omega_{\kappa(d)}$.

Lemma \ref {omegakappa} implies that $\OL\omega_{\kappa(d)}$ is a spherical $3$-design, while Lemma \ref {omegakappa1} implies that $\mathcal D_2(\OL\omega_{\kappa(d)})=\OL\omega_{2d+2}\neq \emptyset$; that is, $\OL\omega_{\kappa(d)}$ is $2$-stiff. Theorem~\ref {2m-1} now implies that the potential $p_f(\ \! \cdot\ \!,\OL\omega_{\kappa(d)})$ attains its absolute minimum on $S^{d-1}$ at every point of $\OL\omega_{2d+2}$.

If the convexity of $f''$ is strict on $(0,4)$, then by Theorem \ref {2m-1}, the set $\mathcal D_2(\OL\omega_{2d+2})=\OL\omega_{\kappa(d)}$ contains all points of absolute minimum of $p_f(\ \!\cdot\ \!,\OL\omega_{2d+2})$ on~$S^{d-1}$, while the set $\mathcal D_2(\OL\omega_{\kappa(d)})=\OL\omega_{2d+2}$ contains all points of absolute minimum of $p_f(\ \!\cdot\ \!,\OL\omega_{\kappa(d)})$ on~$S^{d-1}$.
\end {proof}

\section {A pair of stiff configurations on $S^5$}\label {stiff}

In this section, we consider the symmetrized Schl\"affi configuration $\OL\omega_{54}$ and the (normalized) root system of the $E_6$ root lattice, $\OL\omega_{72}$. Recall that the Schl\"affi configuration, denoted by $\OL \omega_{27}$, consists of the following $N=27$ points (see, e.g., \cite {Boy1995}): the point $(0,0,0,0,0,1)$, ten points whose sixth coordinate is $-1/2$, one of the first five coordinates is $\pm \sqrt{3}/2$, and the remaining four coordinates equal $0$, and sixteen points whose sixth coordinate is $1/4$ and each of the remaining five coordinates is $\pm \sqrt {3}/4$ with an even number of minus signs. The set $\OL\omega_{27}$ is a spherical $4$-design. The symmetrized Schl\"affi configuration is then defined as $\OL\omega_{54}:=\OL\omega_{27}\cup (-\OL\omega_{27})$, which is a spherical $5$-design. 

Denote by $\OL\omega_{72}$ the configuration consisting of $40$ vectors on $S^5$ with two of the first five coordinates equal $\frac {1}{\sqrt{2}}$ or $-\frac {1}{\sqrt{2}}$ and the remaining four coordinates equal zero and $32$ vectors of the form $\(\pm \frac {1}{2\sqrt{2}},\pm\frac {1}{2\sqrt{2}},\pm\frac {1}{2\sqrt{2}},\pm \frac {1}{2\sqrt{2}},\pm \frac {1}{2\sqrt{2}},\pm \frac {\sqrt{3}}{2\sqrt{2}}\)$ with an odd number of minus signs. This is a spherical $5$-design. To show that both configurations are $3$-stiff, we use the following two auxiliary statements.
\begin {lemma}\label {5472}
$\mathcal D_3(\OL\omega_{54})=\OL\omega_{72}$.
\end {lemma}
\begin {proof}
One can verify directly that $\OL\omega_{72}\subset \mathcal D_3(\OL\omega_{54})$. The set $\OL\omega_{54}$ is a $5$-design; that is, $\OL\omega_{54}$ is $3$-stiff. Let now ${\bf x}=(x_1,\ldots,x_6)$ be any vector in $\mathcal D_3(\OL\omega_{54})$. By Proposition \ref {mtight}, we have $D({\bf x},\OL\omega_{54})=\left\{\pm \sqrt{\frac{{3}}{8}},0\right\}$. If ${\bf u}=(0,0,0,0,0,1)\in \OL\omega_{54}$, then $x_6={\bf x}\cdot{\bf u}=\pm \sqrt{\frac {3}{8}}$ or $0$. For every $1\leq i\leq 5$, choose the vector ${\bf v}_i\in\OL\omega_{54}$ whose $i$-th coordinate equals $\frac {\sqrt{3}}{2}$ and the sixth coordinate equals $-\frac {1}{2}$ (all other coordinates are zero). Then 
\begin {equation}\label {3_8}
{\bf x}\cdot{\bf v}_i=\frac {\sqrt{3}}{2}x_i-\frac {x_6}{2}=\pm \sqrt{\frac {3}{8}} \ \ \text{or}\ \  0. 
\end {equation}
If $x_6=0$, then $x_i=\pm \frac {1}{\sqrt{2}}$ or $0$ for $i=1,\ldots,5$. Since $\left|{\bf x}\right|=1$, exactly two coordinates out of first five equal $\pm\frac {1}{\sqrt{2}}$ and four remaining coordinates are zero; that is, ${\bf x}\in \OL\omega_{72}$. 

Assume now that $x_6=\pm\sqrt{\frac {3}{8}}$. Then, for each $1\leq i\leq 5$, from \eqref {3_8}, we have $\frac {\sqrt{3}}{2}x_i=\pm\frac {3}{2}\sqrt{\frac {3}{8}}$ or $\pm\frac {1}{2}\sqrt{\frac {3}{8}}$. Then $x_i=\pm \frac {1}{\sqrt{8}}$, $i=1,\ldots,5$ (the values $x_i=\pm\frac {3}{\sqrt{8}}$ are impossible, since $\left|x_i\right|\leq 1$). Therefore, ${\bf x}=\(\pm \frac {1}{\sqrt{8}},\pm \frac {1}{\sqrt{8}},\pm \frac {1}{\sqrt{8}},\pm \frac {1}{\sqrt{8}},\pm \frac {1}{\sqrt{8}},\pm \sqrt{\frac {3}{8}}\)$. Let ${\bf w}=\(\pm \frac {\sqrt{3}}{4},\pm \frac {\sqrt{3}}{4},\pm \frac {\sqrt{3}}{4},\pm \frac {\sqrt{3}}{4},\pm \frac {\sqrt{3}}{4},\frac {1}{4}\)$ with an even number of minus signs. Assume to the contrary that ${\bf x}$ has an even number of negative coordinates. Then, for each $i=1,\ldots,6$, the product of the $i$-th coordinates of vectors ${\bf x}$ and ${\bf w}$ will be $\pm \frac {\sqrt{3}}{4\sqrt{8}}$ and the number of negative products will be even. This implies that possible values of dot product ${\bf x}\cdot{\bf w}$ will be $\pm \frac {3\sqrt{3}}{2\sqrt{8}}$ or $\pm \frac {\sqrt{3}}{2\sqrt{8}}$, which are not in $D({\bf x},\OL\omega_{54})$. Then ${\bf x}$ must have an odd number of negative coordinates, which makes it an element of $\OL\omega_{72}$.
\end {proof}

\begin {lemma}\label {7254}
$
\mathcal D_3(\OL\omega_{72})=\OL \omega_{54}.
$
\end {lemma}
\begin {proof}
One can verify directly that $\OL\omega_{54}\subset \mathcal D_3(\OL\omega_{72})$ and that $\OL\omega_{72}$ is a $5$-design. Then $\OL\omega_{72}$ is $3$-stiff. Let ${\bf x}=(x_1,\ldots,x_6)$ be any point in $\mathcal D_3(\OL\omega_{72})$. Assume to the contrary that the set $S:=\{\left|x_i\right|: i=1,\ldots,5\}$ has at least two distinct non-zero elements. Then we have $0<\left|x_k\right|<\left|x_\ell\right|$ for some indices $1\leq k\neq \ell\leq 5$. Let ${\bf u}\in \OL \omega_{72}$ be one of the four vectors with the $k$-th and $\ell$-th coordinates being $\pm \frac {1}{\sqrt{2}}$ and all other coordinates being zero. Then dot product ${\bf x}\cdot {\bf u}$ will have four distinct values $\frac{\pm \left|x_k\right|\pm\left|x_\ell\right|}{\sqrt{2}}$ contradicting the assumption that ${\bf x}$ is in $\mathcal D_3(\OL\omega_{72})$. Thus, $x_1=\ldots=x_5=0$ or all non-zero elements of $S$ are equal. If $x_1=\ldots=x_5=0$, then ${\bf x}=(0,0,0,0,0,\pm 1)\in \OL\omega_{54}$.

Assume that all non-zero elements of $S$ are equal and denote by $a$ their common value. Denote also $b:=\left|x_6\right|$. Let $G$ be the set of $32$ vectors $\(\pm \frac {1}{2\sqrt{2}},\pm \frac {1}{2\sqrt{2}},\pm \frac {1}{2\sqrt{2}},\pm \frac {1}{2\sqrt{2}},\pm \frac {1}{2\sqrt{2}},\pm \frac {\sqrt{3}}{2\sqrt{2}}\)$ with an odd number of minus signs.

Assume to the contrary that $b=0$. Let $1\leq k\leq 5$ be such that $\left|x_i\right|=a$ for exactly $k$ indices $i$ between $1$ and $5$ inclusive. Then $a=1/\sqrt{k}$ and the remaining $6-k$ coordinates of ${\bf x}$ are zero.  Since $\OL\omega_{72}$ is $3$-stiff, by Proposition \ref {mtight}, all dot products of ${\bf x}$ with vectors from $G$ must be in the set $\left\{\pm \frac {\sqrt{3}}{2\sqrt{2}},0\right\}$, in particular, there is only one positive dot product. There is a vector ${\bf u}\in G$ with components corresponding to non-zero $x_i$ having the same sign as $x_i$. Then ${\bf x}\cdot {\bf u}=\frac {ka}{2\sqrt{2}}=\frac{\sqrt{k}}{2\sqrt{2}}$, which must equal $\frac {\sqrt{3}}{2\sqrt{2}}$. Consequently, $k=3$. But then there is a vector ${\bf v}\in G$ (with all but one components corresponding to non-zero $x_i$ having the same sign as $x_i$) such that ${\bf x}\cdot {\bf v}=\frac {a}{2\sqrt{2}}=\frac {1}{2\sqrt{6}}$, which is not a possible dot product. Thus, $b$ must be positive.

If $\left|x_i\right|=a$ for exactly one index $1\leq i\leq 5$, then ${\bf x}$, in particular, forms dot products $0$ and $\pm \frac {a}{\sqrt{2}}$ with points from $\OL\omega_{72}$. Since $\OL\omega_{72}$ is $3$-stiff, by Proposition \ref {mtight}, all these dot products must be in the set $\left\{\pm \frac {\sqrt{3}}{2\sqrt{2}},0\right\}$. This implies that $\left|x_i\right|=a=\frac {\sqrt{3}}{2}$ and, since $\left|{\bf x}\right|=1$, we have $x_6=\pm 1/2$. Since all other coordinates of ${\bf x}$ are zero in this case, we have ${\bf x}\in \OL\omega_{54}$.

Assume that $\left|x_i\right|=a$ for exactly $k$ indices $i$ between $1$ and $5$ inclusive, where $2\leq k\leq 4$. Then ${\bf x}$ has at least dot products $0,\pm \frac {a}{\sqrt{2}},\pm \sqrt{2}a$ with vectors from $\OL\omega_{72}$ contradicting the assumption that ${\bf x}\in \mathcal D_3(\OL\omega_{72})$. 

Finally, assume that $\left|x_1\right|=\ldots=\left|x_5\right|=a$. Let ${\bf w}$ be any of the four vectors from $\OL\omega_{72}$ whose first and second coordinates equal $\pm\frac {1}{\sqrt{2}}$ and the remaining four coordinates equal $0$. Then ${\bf x}\cdot {\bf w}$ has values $\pm \sqrt{2}a$ and $0$. Then $\sqrt{2}a=\frac {\sqrt{3}}{2\sqrt{2}}$ or $a=\frac {\sqrt{3}}{4}$. Since $\left|{\bf x}\right|=1$, we have $b=1/4$. 

Assume to the contrary that $x_i<0$ for an odd number of indices $1\leq i\leq 6$. Let $H_1$ be the set of all such indices. For any vector ${\bf u}=(u_1,\ldots,u_6)\in G$, let $H_2:=\{i : u_i<0\}$ (by definition of the set $G$, $\# H_2$ is also odd). For every $i$, we have $\left|x_iu_i\right|=\frac {\sqrt{3}}{8\sqrt{2}}$. Then $x_iu_i=-\frac {\sqrt{3}}{8\sqrt{2}}$ for $\# H_1+\# H_2-2(\#(H_1\cap H_2))$ indices $i$, which is an even number. As ${\bf u}$ runs through the set $G$, dot product ${\bf x}\cdot {\bf u}$ can only assume a value from the set $\left\{\pm \frac {3\sqrt{3}}{4\sqrt{2}},\pm \frac {\sqrt{3}}{4\sqrt{2}}\right\}$, none of which is in the set $\left\{\pm \frac {\sqrt{3}}{2\sqrt{2}},0\right\}$ contradicting Proposition \ref {mtight}. Therefore, coordinates of vector ${\bf x}$ must have an even number of minus signs, which puts vector ${\bf x}$ in $\OL\omega_{54}$. Thus, $\mathcal D_3(\OL\omega_{72})\subset\OL\omega_{54}$.
\end {proof}

We are now ready to prove the main result of this section.
\begin {theorem}\label {72_}
Let $f:[0,4]\to (-\infty,\infty]$ be a function continuous on $(0,4]$ with $f(0)=\lim\limits_{t\to 0^+}f(t)$ and a convex derivative $f^{(4)}$ on $(0,4)$. Then 
\begin {itemize}
\item [(i)]
the potential $p_f(\ \!\cdot\ \!,\OL\omega_{54})$
attains its absolute minimum over $S^{5}$ at every point of $\OL\omega_{72}$;

\item [(ii)]
the potential $p_f(\ \!\cdot\ \!,\OL\omega_{72})$
attains its absolute minimum over $S^{5}$ at every point of~$\OL\omega_{54}$.
\end {itemize}
If the convexity of $f^{(4)}$ is strict on $(0,4)$, then the potential $p_f(\cdot,\OL\omega_{54})$ has no other points of absolute minimum in (i) and the potential $p_f(\ \!\cdot\ \!,\OL\omega_{72})$ has no other points of absolute minimum in (ii).
\end {theorem}
\begin {proof}[Proof of Theorem \ref {72_}]
The set $\OL\omega_{54}$ is a spherical $5$-design. 
In view of Lemma~\ref {5472}, we have $\mathcal D_3(\OL\omega_{54})=\OL\omega_{72}$; that is, $\OL\omega_{54}$ is $3$-stiff. Theorem~\ref {2m-1} now implies that the potential $p_f(\ \! \cdot\ \!,\OL\omega_{54})$ attains its absolute minimum on $S^{5}$ at every point of $\OL\omega_{72}$.
The set $\OL\omega_{72}$ is a spherical $5$-design. Lemma~\ref {7254} implies that $\mathcal D_3(\OL\omega_{72})=\OL\omega_{54}$; that is, $\OL\omega_{72}$ is $3$-stiff. Theorem~\ref {2m-1} now implies that the potential $p_f(\ \! \cdot\ \!,\OL\omega_{72})$ attains its absolute minimum on $S^{5}$ at every point of $\OL\omega_{54}$.

If the convexity of $f^{(4)}$ is strict on $(0,4)$, then by Theorem \ref {2m-1}, the set $\mathcal D_3(\OL\omega_{54})=\OL\omega_{72}$ contains all points of absolute minimum of $p_f(\ \!\cdot\ \!,\OL\omega_{54})$ on~$S^{5}$, while the set $\mathcal D_3(\OL\omega_{72})=\OL\omega_{54}$ contains all points of absolute minimum of $p_f(\ \!\cdot\ \!,\OL\omega_{72})$ on~$S^5$.
\end {proof}

\section {A pair of stiff configurations on $S^6$}\label {S6}

On $S^6$, the kissing arrangement of $N=56$ points, denoted by $\OL\omega_{56}$, was proved in \cite {BoyDraHarSafStosharpantipodal} to be $3$-stiff. We find that its dual is the set of $N=126$ minimal non-zero vectors of $E_7$ root lattice, denoted by $\OL\omega_{126}$, which is also $3$-stiff. The configuration $\OL \omega_{56}$ can be defined as the set of points on the intersection of the sphere $S^7\subset \RR^8$ with the seven-dimensional linear subspace $\mathcal H$ of $\RR^8$ perpendicular to vector $(1,1,1,1,1,1,1,1)$ that have two coordinates equal to $-\frac{3}{2\sqrt{6}}$ and six coordinates equal to $\frac{1}{2\sqrt{6}}$ or two coordinates equal to $\frac {3}{2\sqrt{6}}$ and six coordinates equal to $-\frac{1}{2\sqrt{6}}$. There are only three distinct values of the dot product between any two distinct vectors in $\OL\omega_{56}$, which are $\pm \frac {1}{3}$ and~$-1$. It is a known $3$-sharp configuration. 

The configuration $\OL\omega_{126}$ can be described as the arrangement of $126$ points on the intersection $S^7\cap \mathcal H$ in $\RR^8$ that have one coordinate equal to $\frac {1}{\sqrt {2}}$, one coordinate equal to $-\frac {1}{\sqrt {2}}$, and six zero coordinates or four coordinates equal to $\frac {1}{2\sqrt {2}}$ and four coordinates equal to $-\frac {1}{2\sqrt {2}}$. One can also think of $\OL\omega_{126}$ as the set of minimal non-zero vectors of $E_8$ root lattice whose coordinates sum to zero.

To prove the main result of this section, we use the following two auxiliary statements.
\begin {lemma}\label {56}
$
\mathcal D_{3}(\OL\omega_{56})=\OL\omega_{126}.
$
\end {lemma}
\begin {proof}
It is not difficult to verify directly that $\OL\omega_{126}\subset\mathcal D_3(\OL\omega_{56})$. To show that $\mathcal D_3(\OL\omega_{56})$ is contained in $\OL\omega_{126}$, we will consider these configurations as subsets of the sphere 
$$
\W S^6:=S^7\cap \mathcal H \subset \RR^8.
$$
Let 
$
{\bf x}=(x_1,\ldots,x_8)\in \W S^6$ be arbitrary point in $\mathcal D_3(\OL\omega_{56})$. For every pair of indices $1\leq i\neq j\leq 8$, let ${\bf u}$ be the vector in $\OL\omega_{56}$ with the $i$-th and the $j$-th coordinates equal to $-\frac {3}{2\sqrt {6}}$. Since $x_1+\ldots+x_8=0$, we have
$$
{\bf x}\cdot {\bf u}=\frac {1}{2\sqrt{6}}\(-3x_i-3x_j+\sum_{k:k\neq i,j} x_k\)=-\frac {2}{\sqrt{6}}(x_i+x_j)
$$
and ${\bf x}\cdot(-{\bf u})=\frac {2}{\sqrt{6}}(x_i+x_j)$. Since $\OL\omega_{56}$ is a $5$-design and $\mathcal D_3(\OL\omega_{56})\neq \emptyset$, the set $\OL\omega_{56}$ is $3$-stiff.
Since ${\bf x}\in \mathcal D_3(\OL\omega_{56})$, by Proposition \ref {mtight}, we have 
\begin {equation}\label {V1}
V:=\{\pm (x_i+x_j):1\leq i\neq j\leq 8\}\subset \left\{\pm \frac{1}{\sqrt{2}},0\right\}. 
\end {equation}
Assume that ${\bf x}$ has at least three distinct coordinates, say $x_k<x_\ell<x_m$. Then $V$ contains the numbers $x_k+x_\ell<x_k+x_m<x_\ell+x_m$. This forces $x_k+x_\ell=-\frac {1}{\sqrt{2}}$, $x_k+x_m=0$, and $x_\ell+x_m=\frac {1}{\sqrt{2}}$. The only solution to this system is $x_k=-\frac {1}{\sqrt{2}}$, $x_\ell=0$, and $x_m=\frac {1}{\sqrt{2}}$. If for some $i\neq k,\ell,m$, we had $x_i\neq 0$, then $V$ would contain the fourth element: $x_i+x_m>\frac {1}{\sqrt{2}}$ if $x_i>0$ or $x_i+x_k<-\frac {1}{\sqrt{2}}$ if $x_i<0$. Thus, $x_k=-\frac {1}{\sqrt{2}}$, $x_m=\frac {1}{\sqrt{2}}$, and the remaining six coordinates are $0$; i.e., ${\bf x}\in \OL\omega_{126}$.

Assume now that ${\bf x}$ has at most two distinct coordinates. If all coordinates of ${\bf x}$ were the same, then ${\bf x}$ would have to be a zero vector, since the sum of its coordinates is zero. Thus, ${\bf x}$ has exactly two distinct coordinates. Assume to the contrary that one of the two values of coordinates of ${\bf x}$ appears only once. Denote this value by $c$ and let the other value equal $b$.
Then $c=-7b\neq 0$ and $V$ contains four distinct values: $\pm 6b,\pm 2b$.  This contradicts \eqref {V1}.

Thus, each of the two values of coordinates of ${\bf x}$ appears more than once. Denote them so that $b<c$. Then $V$ contains the numbers $2b<b+c<2c$, which forces $2c=\frac {1}{\sqrt{2}}$ and $2b=-\frac {1}{\sqrt{2}}$; that is, all coordinates of ${\bf x}$ are $\pm\frac {1}{2\sqrt{2}}$. Since the sum of coordinates of ${\bf x}$ is zero, four of its coordinates equal $\frac {1}{2\sqrt{2}}$ and four equal $-\frac {1}{2\sqrt {2}}$; i.e., ${\bf x}\in \OL\omega_{126}$.
\end {proof}

\begin {lemma}\label {126}
$\mathcal D_3(\OL\omega_{126})=\OL\omega_{56}$.
\end {lemma}
\begin {proof}
It is not difficult to verify directly that $\OL\omega_{56}\subset \mathcal D_3(\OL\omega_{126})$. To show that $\mathcal D_3(\OL\omega_{126})$ is contained in $\OL\omega_{56}$, we will consider these configurations as subsets of the sphere $\W S^6$. Let ${\bf x}\in \W S^6$ be an arbitrary point in $\mathcal D_3(\OL\omega_{126})$. For any $1\leq i\neq j\leq 8$, let ${\bf u}\in \OL\omega_{126}$ be the vector whose $i$-th coordinate is $\frac {1}{\sqrt {2}}$, $j$-th coordinate equals $-\frac {1}{\sqrt{2}}$, and the remaining six coordinates equal zero. Then ${\bf x}\cdot{\bf u}=\frac {x_i-x_j}{\sqrt{2}}$. 

Assume to the contrary that vector ${\bf x}$ has at least three distinct coordinates, which we denote by $x_k<x_\ell<x_m$. Then ${\bf x}\cdot {\bf u}$ has, in particular, four distinct values $\pm\frac {x_\ell-x_k}{\sqrt{2}}$ and $\pm \frac{x_m-x_k}{\sqrt{2}}$ contradicting the fact that ${\bf x}\in \mathcal D_3(\OL\omega_{126})$. Thus, there are at most two distinct numbers among coordinates of ${\bf x}.$ If all coordinates of ${\bf x}$ were equal, since they sum to zero, ${\bf x}$ would have to be a zero vector contradicting the assumption that ${\bf x}\in \W S^6$.

Thus, there are exactly two distinct numbers among coordinates of ${\bf x}$, which we denote by $b<c$. By choosing vector ${\bf u}\in \OL\omega_{126}$ with $\frac {1}{\sqrt{2}}$ and $-\frac {1}{\sqrt{2}}$ on appropriate positions and $0$ on the remaining six positions, we get value ${\bf x}\cdot{\bf u}=\frac {c-b}{\sqrt{2}}$. Since ${\bf x}\in \mathcal D_3(\OL\omega_{126})$, by Proposition \ref {mtight}, ${\bf x}\cdot {\bf u}\in \{\pm \frac {1}{\sqrt{3}},0\}$. The only value the positive dot product $\frac {c-b}{\sqrt{2}}$ can have is $\frac {1}{\sqrt{3}}$; that is, $c=b+\sqrt{\frac {2}{3}}$. Denote by $k$ the number of coordinates of ${\bf x}$ that equal $c$. Then $(8-k)b+k\(b+\sqrt{\frac {2}{3}}\)=0$ and $(8-k)b^2+k\(b+\sqrt{\frac {2}{3}}\)^2=1$. From the first equation, we have $b=-\frac {k}{8}\sqrt{\frac {2}{3}}$. Then the second equation implies that $k(8-k)=12$; that is, $k=2$ or $6$. When $k=2$, we have $b=-\frac {1}{4}\sqrt{\frac {2}{3}}=-\frac {1}{2\sqrt{6}}$ and $c=\frac {3}{4}\sqrt{\frac {2}{3}}=\frac {3}{2\sqrt{6}}$. When $k=6$, we have $b=-\frac {3}{4}\sqrt{\frac {2}{3}}=-\frac {3}{2\sqrt{6}}$ and $c=\frac {1}{4}\sqrt{\frac {2}{3}}=\frac {1}{2\sqrt{6}}$. In both cases, ${\bf x}\in \OL\omega_{56}$.
\end {proof}

%
%
%
%

We are now ready to prove the main result of this section.

\begin {theorem}\label {56_}
Let $f:[0,4]\to (-\infty,\infty]$ be a function continuous on $(0,4]$ with $f(0)=\lim\limits_{t\to 0^+}f(t)$ and a convex derivative $f^{(4)}$ on $(0,4)$. Then 
\begin {itemize}
\item [(i)]
the potential $p_f(\ \!\cdot\ \!,\OL\omega_{56})$
attains its absolute minimum over $S^{6}$ at every point of $\OL\omega_{126}$;

\item [(ii)]
the potential $p_f(\ \!\cdot\ \!,\OL\omega_{126})$
attains its absolute minimum over $S^6$ at every point of $\OL \omega_{56}$.
\end {itemize}
If the convexity of $f^{(4)}$ is strict on $(0,4)$, then the potential $p_f(\cdot,\OL\omega_{56})$ has no other points of absolute minimum in (i) and the potential $p_f(\cdot,\OL\omega_{126})$ has no other points of absolute minimum in (ii).
\end {theorem}
We remark that one universal minimum of $\OL\omega_{56}$ and the value of the absolute minimum of its potential on $S^6$ were found in \cite {BoyDraHarSafStosharpantipodal}.

\begin {proof}[Proof of Theorem \ref {56_}]
Since $\OL\omega_{56}$ is a spherical $5$-design and, by Lemma \ref {56}, $\mathcal D_3(\OL\omega_{56})=\OL\omega_{126}\neq \emptyset$, the configuration $\OL\omega_{56}$ is $3$-stiff (see also \cite {BoyDraHarSafStosharpantipodal}). Then we apply Theorem \ref {2m-1} to establish item (i). 
The set $\OL\omega_{126}$ is a spherical $5$-design. Lemma~\ref {126} implies that $\mathcal D_3(\OL\omega_{126})=\OL\omega_{56}\neq \emptyset$; that is, $\OL\omega_{126}$ is $3$-stiff. Theorem~\ref {2m-1} now implies that the potential $p_f(\ \! \cdot\ \!,\OL\omega_{126})$ attains its absolute minimum on $S^{6}$ at every point of $\OL\omega_{56}$.

If, in addition, $f^{(4)}$ is strictly convex on $(0,4)$, Theorem \ref {2m-1} implies that the potential $p_f(\cdot,\OL\omega_{56})$ attains its absolute minimum over $S^6$ only at points of $\OL\omega_{126}$ and the set $\OL\omega_{56}$ contains all points of absolute minimum of $p_f(\ \!\cdot\ \!,\OL\omega_{126})$ on~$S^{6}$.
\end {proof}

}

\begin {thebibliography}{99}
\bibitem {BanDam1979}
E. Bannai, R.M. Damerell, Tight spherical designs I, {\it J. Math. Soc. Japan} {\bf 31} (1979), no. 1, 199--207.
\bibitem {Bil2015}
M. Bilogliadov, Equilibria of Riesz potentials generated by point charges at the roots of unity, {\it Comput. Methods Funct. Theory} {\bf 15} (2015), no. 4, 471--491.
\bibitem {Borsymmetric}
S.V. Borodachov, Absolute minima of potentials of certain regular spherical configurations (submitted), https://arxiv.org/abs/2210.04295.
\bibitem {BorMinMax}
S.V. Borodachov, Min-max polarization for certain classes of sharp configurations on the sphere (submitted), https://arxiv.org/abs/2203.13756.
\bibitem {BorAIMtalk}
S.V. Borodachov, The current state of the minimal dicrete polarization theory, {\it Workshop ``Minimal Energy Problems with Riesz Potentials"}, American Institute of Mathematics, Zoom, May 3--7, 2021.
\bibitem{Bor2022talk}
S.V. Borodachov, Min-max polarization for certain classes of sharp configurations on the sphere, {\it Workshop "Optimal Point Configurations on Manifolds"}, ESI, Vienna, January 17--21, 2022. https://www.youtube.com/watch?v=L-szPTFMsX8
\bibitem {BorBatumi}
S.V. Borodachov, Extreme values of potentials of spherical designs and the polarization problem, {\it XII Annual International Conference of the Georgian Mathematical Union}, Batumi State University, Georgia, August 29--September 3.
\bibitem {Borsimplex}
S.V. Borodachov, Polarization problem on a higher-dimensional sphere for a simplex, {\it Discrete and Computational Geometry} {\bf 67} (2022), 525--542.
\bibitem{BorHarSafbook}
S. Borodachov, D. Hardin, E. Saff, {\it Discrete Energy on Rectifiable Sets}. Springer, 2019.
\bibitem {BoyDraHarSafSto600cell}
P.G. Boyvalenkov, P.D. Dragnev, D.P. Hardin, E.B. Saff, M.M. Stoyanova, On polarization of spherical codes and designs (submitted), https://arxiv.org/abs/2207.08807.
\bibitem {BoyDraHarSafStosharpantipodal}
P.G. Boyvalenkov, P.D. Dragnev, D.P. Hardin, E.B. Saff, M.M. Stoyanova, Universal minima of discrete potentials for sharp spherical codes, https://arxiv.org/pdf/2211.00092.pdf
\bibitem {Boy1995}
P. Boyvalenkov, 
Computing distance distributions of spherical designs,
{\it Linear Algebra Appl.} {\bf 226/228} (1995), 277--286.
\bibitem {BoyDanKaz2001}
P. Boyvalenkov, D. Danev, P. Kazakov, Indexes of spherical codes, {\it Codes and association schemes} (Piscataway, NJ, 1999), 47–57. {\it DIMACS Ser. Discrete Math. Theoret. Comput. Sci.}, {\bf 56}, Amer. Math. Soc., Providence, RI, 2001.
\bibitem{CohKum2007}
H. Cohn, A. Kumar, Universally optimal distribution of points on spheres, {\it J. Amer. Math. Soc.} {\bf 20} (2007), no. 1, 99--148.
\bibitem {Dav1975}
P.J. Davis, {\it Interpolation and Approximation}, Second edition, Dover Publications, New York, NY, 1975. 
\bibitem {DelGoeSei1977}
P. Delsarte, J.M. Goethals, J.J. Seidel, Spherical codes and designs, {\it Geometriae Dedicata}, {\bf 6} (1977), no. 3, 363--388.
\bibitem {Epp2013}
J.F. Epperson, {\it An introduction to numerical methods and analysis},  Wiley, 2013.
\bibitem {GioKhi2020}
G. Giorgadze, G. Khimshiashvili, Stable equilibria of three constrained unit charges, {\it Proc. I. Vekua Inst. Appl. Math.} {\bf 70} (2020), 25--31.
\bibitem {HarKenSaf2013}
D. Hardin, A. Kendall, E. Saff,
Polarization optimality of equally spaced points on the circle for discrete potentials.
{\it Discrete Comput. Geom.} {\bf 50} (2013), no. 1, 236--243. 
\bibitem {IsaKel1965}
E. Isaacson, H. Keller, {\it Analysis of numerical methods}. Dover Books, 1994.
\bibitem {Lev1979}
V.I.Levenshtein, On bounds for packings in $n$-dimensional Euclidean space,
{\it Soviet Math. Dokladi}, 20, 1979, 417--421.
\bibitem {Lev1992}
V.I. Levenshtein, Designs as maximum codes in polynomial metric spaces,
{\it Acta Appl. Math.} {\bf 25} (1992), 1--82.
\bibitem {Lev1998}
V.I. Levenshtein, Universal bounds for codes and designs, Chapter 6 in
{\it Handbook of Coding Theory}, V. Pless and W.C. Huffman, Eds., Elsevier Science B.V., 1998.
\bibitem {Mus2003}
O.R. Musin, The kissing number in four dimensions, {\it Annals of Mathematics}, {\bf 168} (2008), 1--32.
\bibitem {NikRaf2011}
N. Nikolov, R. Rafailov, On the sum of powered distances to certain sets of points on the circle, {\it Pacific J. Math.} {\bf 253} (2011), no. 1, 157--168.
\bibitem{NikRaf2013}
N. Nikolov, R. Rafailov, On extremums of sums of powered distances to a finite set of points. {\it Geom. Dedicata} {\bf 167} (2013), 69--89.
\bibitem {NIST}
{\it NIST Digital Library of Mathematical Functions.} http://dlmf.nist.gov/,
Release 1.0.13 of 2016-09-16. F. W. J. Olver, A. B. Olde Daalhuis, D. W.
Lozier, B. I. Schneider, R. F. Boisvert, C. W. Clark, B. R. Miller and
B. V. Saunders, eds.
\bibitem {Sto1975circle}
K. Stolarsky, The sum of the distances to certain pointsets on the unit circle, {\it Pacific J. Math.} {\bf 59} (1975), no. 1, 241--251. 
\bibitem {Sto1975}
K. Stolarsky, The sum of the distances to $N$ points on a sphere, {\it Pacific J. Math.} {\bf 57} (1975), no. 2, 563--573.
\bibitem {Sze1975}
G. Szeg\"o, {\it Orthogonal polynomials}. Fourth edition. American Mathematical Society, Colloquium Publications, Vol. XXIII. American Mathematical Society, Providence, R.I., 1975.
\end {thebibliography}
\end {document}